\documentclass{article} 

\usepackage{amsfonts,amsmath,amssymb}
\usepackage[all]{xy} 
\usepackage{epsfig, latexsym}

\usepackage[latin1]{inputenc}

\newenvironment{preuve}[1][]
{\vskip 2mm\noindent {\bf D\'emonstration #1:  }}{$\Box$ \vskip 2mm}
\newenvironment{theor}
{\vskip 2mm\noindent {\bf Th\'eor\`eme.}}{\vskip 2mm}
\newenvironment{defin}
{\vskip 2mm\noindent {\bf D\'efinition.}}{\vskip 2mm}

\newtheorem{theorem}{Th\'eor\`eme}
\newtheorem{lemma}{Lemme}
\newtheorem{corollaire}{Corollaire}
\newtheorem{definition}{D\'efinition}
\newtheorem{prop}{Proposition}
\newtheorem{proposition-definition}{Proposition-D\'efinition}
\def\min{\mathop{\mathrm{inf}}}

\newcommand{\calv}{\mathcal{V}}

\newcommand{\calf}{\mathcal{F}}

\newcommand{\calu}{\mathcal{U}}
\newcommand{\cale}{\mathcal{E}}

\title{Niveaux de  ${\bf R}$-scindement de structures de Hodge mixtes }

\begin{document}

\author {\normalsize{Olivier Penacchio }}

\maketitle
\centerline{\textbf{Abstract}}
We propose a generalization of the notion of ${\bf R}$-split mixed Hodge structure by defining a ${\bf R}$-splitting level for mixed Hodge structures. This is a discrete invariant taking values in positive integers and equal to $0$ for ${\bf R}$-split mixed Hodge structures. We describe some properties of this invariant for classical operations on these structures and give some explicit calculation.\\
                                                                                                                                                                                                                                                                                                                                                                                                                     
\centerline{\textbf{R\'esum\'e}}
Nous proposons une généralisation de la notion de structure de Hodge ${\bf R}$-scindée en définissant un niveau de ${\bf R}$-scindement pour les structures de Hodge mixtes. C'est un invariant discret à valeur dans les entiers positifs qui vaut $0$ pour une structure ${\bf R}$-scindée. Nous décrivons le comportement de cet invariant par les opérations naturelles sur les structures de Hodge mixtes et donnons quelques exemples de calculs. \\

\textsc{Mots-clefs}: Théorie de Hodge, Matrice des périodes, Fibrés vectoriels.

\textsc{Code mati\`ere AMS}: 14C30-32G20.

 \section{Introduction}

Les structures de Hodge mixtes apparaissent naturellement sur les groupes de cohomologie des variétés complexes (voir \cite{del2},\cite{del3}). Ce sont, pour ainsi dire, des extensions successives de structures de Hodge pures, mais elles n'admettent pas en général de décomposition en somme directe par des sous-espaces $H^{p,q}$ tels que $H^{p,q}={\overline H^{q,p}}$. On peut cependant trouver une décomposition qui a des propriétés analogues. En effet, d'après Deligne, la donnée d'une structure de Hodge mixte $(W_{\bullet},F^{\bullet })$ sur un espace vectoriel $H_{\bf C}=H_{\bf R}\otimes_{\bf R}{\bf C}$ définit une bigraduation canonique de $H_{\bf C}$ : $H_{\bf C}=\oplus_{p,q}\,I^{p,q}$ qui satisfait les propriétés suivantes :\\

\hspace{1.5cm}$(i)$ $W_{m}=\oplus_{p+q \leq m}I^{p,q}$.

\hspace{1.5cm}$(ii)$ $F^{p}=\oplus_{i \geq p} \oplus_{q} I^{i,q}$.
 
\hspace{1.5cm}$(iii)$ $I^{p,q}={\overline I}^{q,p} \ \ \text{mod} \ W_{p+q-2}. $\\
\\
La structure de Hodge mixte est dite ${\bf R}$-scindée si pour tous $(p,q)$ : $I^{p,q}=\overline{I^{q,p}}$. C'est le cas en particulier, d'après $(iii)$, si c'est une structure de Hodge pure ou si la filtration par le poids $W_{\bullet }$ est de longueur $1$.

D'après \cite{catkapsch},  à toute structure de Hodge mixte peut \^etre associée de manière canonique une structure de Hodge mixte ${\bf R}$-scindée en modifiant la filtration de Hodge $F^{\bullet }$. Ce résultat se montre en considérant l'algèbre de Lie nilpotente :\hspace{0.5cm}${\mathfrak l}^{-1,-1}=\{X \in \text{End}(V_{\bf C}): \ X(I^{p,q}) \subset \oplus_{r<p,s<q}I^{r,s}\}$\\
Le fait que ${\overline {\mathfrak l}^{-1,-1}}={\mathfrak l}^{-1,-1}$ induit une structure réelle sous-jacente, ${\mathfrak l}^{-1,-1}_{\bf R}=\text{End}(V_{\bf R})\cap {\mathfrak l}^{-1,-1}$. On montre alors que pour toute structure de Hodge mixte $H=(H_{\bf C},W_{\bullet },F^{\bullet })$ il existe une unique $\delta \in {\mathfrak l}^{-1,-1}_{\bf R}$ tel que $H'=(H_{\bf C},W_{\bullet },e^{-i.\delta}.F^{\bullet })$ soit ${\bf R}$-scindée.

Notre intention ici est de mesurer, par un invariant à valeurs discretes, à quel point une structure de Hodge mixte $H=(H_{\bf C},W_{\bullet },F^{\bullet })$ est "loin" de la structure ${\bf R}$-scindée associée en utilisant les positions relatives des trois filtrations définissant la structure de Hodge mixte : la filtration par le poids, la filtration de Hodge et la filtration conjuguée à la filtration de Hodge par rapport à la structure réelle sous-jacente. Dans ce but, on considère les dimensions des intersections des sous-espaces donnés par les filtrations. Tout d'abord, les nombres de Hodge $h^{p,q}_{H}=\text{dim}_{\bf C}Gr^{F^{\bullet }}_{p}Gr_{W_{\bullet }}^{p+q}H_{\bf C}$. Ils sont constants sur l'ensemble des structures de Hodge mixtes obtenues par action sur $H$ de l'algèbre de Lie nilpotente ${\mathfrak l}^{-1,-1}$. Nous introduisons les entiers $t^{p,q}_{H}=\text{dim}_{\bf C}Gr^{{\overline F}^{\bullet }}_{q}Gr^{F^{\bullet }}_{p}H_{\bf C}$. En général ces entiers ne sont pas constants lors de la déformation de $H$ par l'action considérée. Leur définition fait intervenir la filtration conjuguée à la filtration de Hodge, ainsi, pour une variation de structures de Hodge mixte, ils sont constants sur des strates analytiques réelles de la base de la variation (voir \cite{pen}, en préparation). La motivation géométrique de l'introduction d'un tel invariant est de permettre d'étudier les variations de structures de Hodge mixtes à l'aide de ces strates analytiques réelles. On donne ainsi comme définition du niveau de ${\bf R}$-scindement : 
\begin{defin} Le niveau de ${\bf R}$-scindement d'une structure de Hodge mixte $H=(H_{\bf C},W_{\bullet },F^{\bullet })$ est l'entier :$$\alpha(H)=\frac{1}{2}\sum_{p,q}\ (p+q)^{2}.(h^{p,q}_{H}-t^{p,q}_{H}).$$
\end{defin}
Lorsque $H$ est ${\bf R}$-scindée, comme $t^{p,q}_{H}=h^{p,q}_{H}=\text{dim}_{\bf C}\ I^{p,q}$, l'invariant $\alpha$ est nul. En général, c'est un entier négatif. Avant de décrire le comportement de cet invariant par les opérations naturelles sur les structures de Hodge, expliquons qu'il peut \^etre obtenu à partir du caractère de Chern d'un fibré vectoriel sur ${\bf P}^2$ construit à partir des filtrations. 
 
Considérons le cas général où l'espace filtré que l'on étudie ne provient pas forcément d'une structure de Hodge mixte. Soit, donc, $(V,F^{\bullet }_{0},F^{\bullet }_{1},F^{\bullet }_{2})$ un espace vectoriel de dimension finie muni de trois filtrations décroissantes et complètes. Pour étudier les positions relatives de ces filtrations sur $V$, on associe à ce triplet de filtrations un fibré vectoriel $\xi_{{\bf P}^2}(V,F^{\bullet }_{0},F^{\bullet }_{1},F^{\bullet }_{2})$ sur ${\bf P}^2$, appelé fibré de Rees. Il est obtenu en recollant des faisceaux cohérents sur les cartes affines standards ${\bf A}^{2}\cong\text{Spec}{\bf C}[u,v]$ recouvrant le plan projectif associés à des ${\bf C}[u,v]$-modules construits à partir des intersections des sous-espaces associés aux filtrations. On peut alors de explicitement calculer le caractère de Chern $\text{ch}(\xi_{{\bf P}^2}(V,F^{\bullet }_{0},F^{\bullet }_{1},F^{\bullet }_{2}))$ à partir des dimensions des intersections entre ces sous-espaces. Revenons aux structures de Hodge mixtes. Soit $H=(H_{\bf C},W_{\bullet },F^{\bullet })$ une telle structure, on applique la construction du fibré de Rees à $(H_{\bf C},W^{\bullet }_{.},F^{\bullet },{\overline F}^{\bullet})$, où $W^{\bullet }_{.}$ est la filtration décroissante associée à la filtration croissante $W_{\bullet }$. On peut alors montrer que $\alpha(H)=\text{c}_{2}(H_{\bf C},W^{\bullet }_{.},F^{\bullet },{\overline F}^{\bullet })$. La motivation pour exprimer l'invariant ainsi est de pouvoir étudier son comportement en famille en associant une famille de fibrés vectoriels à une variation de structures de Hodge mixte. On utilise ensuite la semi-continuité de cet invariant pour définir une stratification analytique réelle de la base de la variation ou des espaces classifiants de structures de Hodge mixtes en général (voir \cite{pen}).

L'invariant $\alpha$ se comporte bien vis-à-vis des opérations naturelles sur les structures de Hodge mixtes : 
\begin{theor}
Soient $H$ et $H'$ deux structures de Hodge mixtes, alors :

$(i)$ pour tout $k \in {\bf Z}$, $\alpha(H \otimes T \langle k \rangle )=\alpha(H)$. 

$(ii)$ $\alpha(H^{*})=\alpha(H)$ où $H^{*}=\emph{Hom}_{SHM}(H,T \langle 0 \rangle )$.

$(iii)$ $\alpha ( H \oplus H')=\alpha(H) +\alpha (H')$.
 
$(iv)$ $\alpha( H \otimes H')=\emph{dim}(H').\alpha(H)+\emph{dim}(H).\alpha(H')$.\\ 
(où $T \langle k \rangle$ est la structure de Hodge de Tate de type $(-k,-k)$).
\end{theor}

Il est intéressant géométriquement d'avoir des informations sur le comportement de l'invariant par des morphismes des structures de Hodge. On montre par exemple que si $f: X \rightarrow Y$ est un morphisme entre variétés algébriques complexes qui induit un morphisme injectif sur la cohomologie, alors pour tout $k$, $\alpha(H^{k}(X,{\bf C})) \geq \alpha(H^{k}(Y,{\bf C}))$. Notons $\text{Ext}(B,A)$ le groupe des classes de congruence d'extensions de structures de Hodge mixtes $H$ de $B$ par $A$. L'assertion précédente est une conséquence directe de la sur-additivité de l'invariant par extension :

\begin{theor}\label{inegext}
Soient $A$ et $B$ deux structures de Hodge mixtes et $H \in \emph{Ext}(B,A)$, alors :
$$ \alpha(H) \geq \alpha(A)+\alpha(B)$$

\end{theor} 
Pour finir, on donne le calcul explicite de de $\alpha(.)$ pour les structures de Hodge mixtes sur les groupes de cohomologie des courbes singulières et non complètes en genre $0$ et $1$ en déduisant les entiers $t^{p,q}_{H}$ de la matrice des périodes.\\

\textbf{Remerciements.} Je tiens à remercier Carlos Simpson, qui m'a donné ce sujet de recherche, pour ses encouragements et son aide, et Philippe Eyssidieux qui m'a accordé de nombreuses discussions.


\section{Fibré de Rees sur ${\bf P}^{2}$ associé à trois filtrations }

 Dans cette section nous allons associer un fibré vectoriel sur ${\bf P}^{2}$ à tout espace vectoriel muni de trois filtrations $(V,F^{\bullet }_{1},F^{\bullet }_{2},F^{\bullet }_{3})$.  L'étude des invariants de ce fibré permettra d'avoir des informations sur les positions relatives des filtrations dans $V$. 
\subsection{Filtrations et module de Rees}
Soit $V $ un espace vectoriel complexe de dimension finie muni d'une filtration décroissante $F^{\bullet }$ c'est à dire d'une suite de sous-espaces vectoriels indexés par ${\bf Z} $ telle que : $\forall m,n \,\,\, n\leq m \Rightarrow  F^{n} \subset F^{n}$. Une filtration $F^{\bullet }$ sur $V$ est dite complète s'il existe deux entiers $p$ et $q$ tels que $F^{p}=V$ et $F^{q}=\{ 0 \}$.

Tous les espaces vectoriels complexes munis de filtrations seront supposés de dimension finie. Toutes les filtrations seront supposées complètes et sauf avis contraire décroissantes. On note $End(V)$ l'anneau des endomorphisme de $V$.
 
\begin{definition}
Soit $(V,F^{\bullet })$ un espace vectoriel filtré. Un endomorphisme semi-simple $Y \in gl(V)$ scinde $F^{\bullet }$ si $\forall k \in {\bf Z} \,\,\, F^{k}=F^{k+1}\oplus E_{k}(Y)$ où $E_{k}(Y)$ est le sous-espace propre associé à la valeur propre $k$. $Y$ est une graduation de $(V,F^{\bullet })$.
\end{definition}
Toute filtration $F^{\bullet }$ sur un espace vectoriel $V$ peut \^etre scindée. On peut construire un scindement de la façon suivante : $F^{\bullet }$ est décroissante et complète donc il existe un plus grand indice $p$ tel que $F^{p+1}=\{0\}$ et $F^{p} \neq \{0\}$. Soit $\{f^{p}_{i}\}_{i \in I_p}$ une famille d'éléments de $V$ formant une base de $F^p$. On complète cette famille dans $F^{p-1}$ par des vecteurs $\{f^{p-1}_{i}\}_{i \in I_{p-1}}$ pour obtenir une base de $F^{p-1}$. On itére l'opération pour obtenir une base de $V$ formée par la famille finie $\cup_{k \in [q,p]}  \{f^{k}_{i}\}_{i \in I_k}$ où $q$ est un entier tel que $F^{q}=V$. Cette famille est dite compatible avec la filtration. Soit $Y$ l'endomorphisme de $V$ définit par $Y(f^{k}_{i})=k.f^{k}_{i}$ pour $i \in I_{k}$. $Y$ est bien un scindement de $(V,F^{\bullet })$. 

Une tel scindement définit une graduation de $V$ i.e une décomposition de $V$ en une somme directe compatible avec la filtration : $V=\oplus_{k=q}^{k=p}E_{k}(Y)$ et $F^{n}=\oplus_{k=n}^{k=p}E_{k}(Y)$. Nous noterons ${\mathcal Y}(F^{\bullet })$ l'ensemble des endomorphismes de $V$ qui scindent $F^{\bullet }$. \\

Soit $(V,F^{\bullet }_{1},F^{\bullet }_{2})$ un espace vectoriel muni de deux filtrations. Une bigraduation de $V$, $V=\oplus_{p,q}V^{p,q}$, est dite compatible aux filtrations $F^{\bullet }_{1}$ et $F^{\bullet }_{2}$ si pour tout $p \in {\bf Z}$ :

\begin{center}
 $F^{p}_{1}=\oplus_{a \geq p,q}V^{a,q}$ et $F^{p}_{2}=\oplus_{ p,q \geq b}V^{p,b}$. 
\end{center}
Il existe toujours une graduation de $V$ compatible à $F^{\bullet }_{i}$ associée à un endomorphisme $Y_{F^{\bullet }_{i}}$ qui scinde $F^{\bullet }_{i}$ pour $i \in \{1,2\}$. Le lemme suivant signifie que l'on peut trouver deux tels endomorphismes $Y_{F^{\bullet }_{1}} \in {\mathcal Y}(F^{\bullet }_{1})$ et $Y_{F^{\bullet }_{2}} \in {\mathcal Y}(F^{\bullet }_{2})$ qui commutent.

\begin{lemma}\label{bigraduation}
Pour tout espace vectoriel muni de deux filtrations $(V,F^{\bullet }_{1},F^{\bullet }_{2})$ il existe toujours une bigraduation $V=\oplus_{p,q}V^{p,q}$ compatible aux deux filtrations.  
\end{lemma}

On dit que l'on peut toujours scinder deux filtrations simultanément. La preuve consiste à exhiber une base de $V$ compatible avec les deux filtrations par complétions de bases successives. 

\begin{preuve}
Les deux filtrations étant complètes, on peut trouver deux entiers $p$ et $q$ tels que $F^{q }_{1}=F^{q }_{2})=V$ et $F^{p}_{1}=F^{p }_{2}=\{0\}$. Soit $\{f^{k,l}_{i}\}_{i \in I_{p-1,p-1}}$ une famille de vecteurs de $V$ qui forme une base de $F^{p-1}_{1}\cap F^{p-1}_{2}$. On la complète en une base de $F^{p-1 }_{1}\cap F^{p-2 }_{2}$ par une famille $\{f^{k,l}_{i}\}_{i \in I_{p-1,p-2}}$. On peut procéder ainsi jusqu'à $F^{p-1}_{1}\cap F^{q}_{2}=F^{p-1}$. Complétons la base de $F^{p-1}_{1}$ à une base de $F^{p-1}_{1}\cup (F^{p-2}_{1} \cap F^{p-2})$ par une famille $\{f^{k,l}_{i}\}_{i \in I_{p-2,p-1}}$. Les indices $(k,l) \in [q,p-1]\times[q,p-1]$ peuvent ainsi \^etre "descendus" de $(p-1,p-1)$ à $(q,q)$ par l'algorithme suivant : tant que $l>q$ on complète la base de $F^{k+1}_{1}\cup (F^{k}_{1}\cap F^{l}_{2})$ en une base de $F^{k+1}_{1}\cup (F^{k}_{1}\cap F^{l-1}_{2})$ par une famille de vecteurs $\{f^{k,l}_{i}\}_{i \in I_{k,l-1}}$. Si $l=q$, on complète la base de $F^{k}_{1}$ en une base de $F^{k}_{1}\cup (F^{k-1}_{1}\cap F^{p-1}_{2})$ par une famille de vecteurs $\{f^{k,l}_{i}\}_{i \in I_{k-1,p-1}}$. Soient les endomorphismes $Y_{F^{\bullet }_{1}} \in {\mathcal Y}(F^{\bullet }_{1})$ et $Y_{F^{\bullet }_{2}} \in {\mathcal Y}(F^{\bullet }_{2})$ définis pour tout $i \in I_{k,l}$ par  $Y_{F^{\bullet }_{1}}(f^{k,l}_{i})=k.f^{k,l}_{i}$ et $Y_{F^{\bullet }_{1}}(f^{k,l}_{i})=l.f^{k,l}_{i}$. Ces deux endomorphismes commutent et scindent bien les deux filtrations. 
\end{preuve}
 Une telle propriété n'est pas vraie en général lorsque $V$ est équipé de plus de deux filtrations. Soit $(V,F^{\bullet}_{1},...,F^{\bullet }_{n})$ un espace vectoriel muni de $n$ filtrations. Il n'est pas toujours possible de trouver $Y_{F^{\bullet }_{i}} \in {\mathcal Y}(F^{\bullet }_{i})$ pour tous $i \in [1,n]$ tels que pour tous $(i,j) \in [1,n]^{2}$, $Y_{F^{\bullet }_{i}}.Y_{F^{\bullet }_{j}}=Y_{F^{\bullet }_{j}}.Y_{F^{\bullet }_{i}}$. 

\textit{Exemple :}
Soit $V={\bf C}^2=<e,f>_{\bf C}$ équipé des filtrations décroissantes et complètes $F^{\bullet}_{1}$, $F^{\bullet}_{2}$ et $F^{\bullet}_{3}$ suivantes :\\ 
 \hspace*{.5cm}
$\bullet \ F^{0}_{1}=V, \ F^{1}_{1}=<f+\kappa e>, \  F^{2}_{1}=\{0\} $ où $\kappa \in {\bf C}$,\\
                \hspace*{.5cm}
$\bullet \ F^{0}_{2}=V,\  F^{1}_{2}=<f+\lambda e>,\  F^{2}_{2}=\{0\}$ où $\lambda \in {\bf C}$,\\
\hspace*{.5cm}
$\bullet \  F^{0}_{3}=V,\  F^{1}_{3}=<f+\mu e>,\  F^{2}_{3}=\{0\}$ où $\mu \in {\bf C}$.\\
Supposons que $\kappa, \lambda , \mu$ soient distincts deux à deux, sans quoi on est ramené à la situation de scinder au plus deux filtrations. Scinder $F^{\bullet }_{1}$ et $F^{\bullet }_{2}$ simultanément revient à prendre pour bigraduation de $V=V^{1,0} \oplus V^{0,1}$ où $V^{1,0}=<f+\kappa e>$ et $V^{0,1}=<f+\lambda e>$. $V^{1,0}$ est le sous-espace propre associé à la valeur $1$ de tout élément dans ${\mathcal Y}(F^{\bullet }_{1})$, ici on doit choisir (pour que les endomorphismes commutent) pour $Y_{F^{\bullet }_{1}}$ l'élément qui a pour sous-espace propre associé à la valeur $0$, $V^{0,1}$. De m\^eme, $V^{0,1}$ est le sous-espace propre associé à la valeur $1$ de tout élément dans ${\mathcal Y}(F^{\bullet }_{2})$, ici on doit choisir pour $Y_{F^{\bullet }_{2}}$ l'élément qui a pour sous-espace propre associé à la valeur $0$, $V^{1,0}$. Ainsi $Y_{F^{\bullet }_{1}}.Y_{F^{\bullet }_{2}}-Y_{F^{\bullet }_{2}}.Y_{F^{\bullet }_{1}}=0$. Comme $\mu \neq \kappa$ et $\mu \neq \lambda$, on ne peut diagonaliser aucun élément $Y_{F^{\bullet }_{3}} \in {\mathcal Y}(^{\bullet }_{3})$ dans cette base.\\

A partir d'un espace vectoriel muni de $n$ filtrations, nous allons associer un module appelé module de Rees d'ordre $n$
construit à partir des intersections des sous-espaces donnés par les $n$-filtrations.

\begin{definition}
Soit $(V,F^{\bullet}_{1},F^{\bullet}_{2},...,F^{\bullet}_{n})$ un espace vectoriel muni de $n$ filtrations. Le module de Rees d'ordre $n$ associé à $(V,F^{\bullet}_{1},F^{\bullet}_{2},...,F^{\bullet}_{n})$, noté $R^{n}(V,F^{\bullet}_{i})$, est le sous-${\bf C}[u_{1},u_{2},...,u_{n}]$-module de \\
${\bf C}[u_{1},u_{2},...,u_{n},u_{1}^{-1},u_{2}^{-1},...,u_{n}^{-1}]\otimes V$ suivant :
$$R^{n}(V,F^{\bullet}_{i})=\sum_{(p_{1},p_{2},...,p_{n}) \in {\bf Z}^{n}} \,u_{1}^{-p_{1}}u_{2}^{-p_{2}}...\,u_{n}^{-p_{n}}\,(F^{p_{1}}_{1} \cap F^{p_{2}}_{2} \cap ... \cap F^{p_{n}}_{n}).$$
\end{definition} 

Soit $(V,F^{\bullet })$ un espace vectoriel filtré. La filtration décalée de $r$ de $F^{\bullet }$ sur $V$ est la filtration notée $Dec^{r}F^{\bullet }$ donnée pour tout $p \in {\bf Z}$ par $Dec^{r}F^{p}=F^{p-r}$. La filtration triviale sur $V$ est notée $Triv^{\bullet }$ et est définie par $Triv^{0}=V$ et $Triv^{1}=\{0\}$.

On vérifie aisément que :

\begin{lemma}\label{isosummodrees}
Soient $(V,F^{\bullet}_{1},F^{\bullet}_{2},...,F^{\bullet}_{n}) $, $(V',{F^{\bullet}_{1}}',{F^{\bullet}_{2}}',...,{F^{\bullet}_{n}}') $ deux espaces vectoriels munis de $n$ filtrations, alors :
\noindent
\item{(i)} Soit $(V \oplus V',F^{\bullet}_{1}\oplus {F^{\bullet}_{1}}',F^{\bullet}_{2}\oplus{ F^{\bullet}_{2}}',...,F^{\bullet}_{n}\oplus {F^{\bullet}_{n}}')$ l'objet qui s'en déduit par somme directe, on a l'isomorphisme de ${\bf C}[u_{1},u_{2},...,u_{n}]$-modules $ R^{n}(V,F^{\bullet}_{i}) \oplus R^{n}(V',{F^{\bullet}_{i}}') \cong R^{n}(V\oplus V',F^{\bullet}_{i}\oplus {F^{\bullet}_{i}}')$.
\item{(ii)} Soit $(V \otimes V',F^{\bullet}_{1}\otimes {F^{\bullet}_{1}}',F^{\bullet}_{2}\otimes {F^{\bullet}_{2}}',...,F^{\bullet}_{n}\otimes {F^{\bullet}_{n}}')$ l'objet qui s'en déduit par produit tensoriel, on a l'isomorphisme de ${\bf C}[u_{1},u_{2},...,u_{n}]$-modules $ R^{n}(V,F^{\bullet}_{i}) \otimes R^{n}(V',{F^{\bullet}_{i}}') \cong R^{n}(V\otimes V',F^{\bullet}_{i}\otimes {F^{\bullet}_{i}}')$.
\item{(iii)} Soit $(p_{i}) \in {\bf Z}^{n}$, alors : $(\prod_{i \in [1,n]}u_{i}^{p_i}) R^{n}(V,F^{\bullet}_{i}) \cong R^{n}(V,Dec^{-p_{i}}F^{\bullet}_{i})$ est un isomorphisme de ${\bf C}[u_{1},u_{2},...,u_{n}]$-modules.

\end{lemma}

\subsection{Faisceau cohérent sur l'espace affine ${\bf A}^{n}$ associé à un module de Rees d'ordre $n$}

Soit $M$ un $A$-module où $A$ est un anneau commutatif unitaire. Nous allons appliquer la construction d'un faisceau cohérent $\tilde M$ sur $Spec \ A$ au ${\bf C}[u_{1},u_{2},...,u_{n}]$-module $R^{n}(V,F^{\bullet}_{i})$. 

\begin{definition}
Le faisceau cohérent de Rees associé à un espace vectoriel muni de $n$ filtrations $(V,F^{\bullet}_{1},F^{\bullet}_{2},...,F^{\bullet}_{n}) $ est le faisceau cohérent sur ${\bf A}^{n}=Spec \ {\bf C}[u_{1},u_{2},...,u_{n}]$ associé au ${\bf C}[u_{1},u_{2},...,u_{n}]$-module $R^{n}(V,F^{\bullet}_{i})$. Il est noté $\xi_{{\bf A}^n}(V,F^{\bullet}_{i})$.
\end{definition} 

\begin{lemma}
$\xi_{{\bf A}^n}(V,F^{\bullet}_{i})$ est un faisceau cohérent reflexif.
\end{lemma}

Comme l'ensemble singulier d'un faisceau reflexif est de codimension au moins trois (cf \cite{oss} par exemple), tous les faisceaux cohérents reflexifs sur des surfaces sont localement libres. Ainsi, pour tout espace vectoriel muni de deux filtrations $(V,F^{\bullet}_{1},F^{\bullet}_{2}) $, $\xi_{{\bf A}^2}(V,F^{\bullet}_{1},F^{\bullet }_{2})$ est un faisceau localement libre sur ${\bf A}^{2}$ c'est à dire, comme sur une variété algébrique il y a équivalence entre les classes d'isomorphisme de fibrés localement libre et classes d'isomorphisme de fibrés vectoriels, $\xi_{{\bf A}^2}(V,F^{\bullet}_{1},F^{\bullet }_{2})$ est un fibré vectoriel sur ${\bf A}^{2}$. C'est le fibré de Rees sur ${\bf A}^{2}$ associé à $(V,F^{\bullet}_{1},F^{\bullet}_{2}) $.

On peut voir directement, sans utiliser le fait que $\xi_{{\bf A}^2}(V,F^{\bullet}_{1},F^{\bullet }_{2})$ est reflexif, que c'est un fibré vectoriel sur ${\bf A}^2$. Il faut pour cela utiliser la bigraduation  de $V$ compatible aux deux filtrations exhibée au lemme \ref{bigraduation}. Soit $V=\oplus_{p,q}V^{p,q}$ une bigraduation compatible aux filtrations. On a une description directe du module de Rees $R^{2}(V,F^{\bullet}_{1},F^{\bullet}_{2})$ en fonction de cette bigraduation, c'est le sous ${\bf C}[u_{1},u_{2}]$-module de ${\bf C}[u_{1},u_{2},u_{1}^{-1},u_{2}^{-1}] \otimes V$ suivant :
\begin{center}
$R^{2}(V,F^{\bullet}_{1},F^{\bullet }_{2})=\oplus_{(p,q) \in {\bf Z}^{2}} \,u_{1}^{-p}u_{2}^{-q}\,V^{p,q},$ 
\end{center}    
où la somme est finie. Soit $\phi$ le morphisme de ${\bf C}[u_{1},u_{2}]$-modules de $R^{2}(V,F^{\bullet}_{1},F^{\bullet }_{2})$ vers $\oplus_{(p,q) \in {\bf Z}^{2}}V^{p,q}$ donné par $\phi(v^{p,q})=u_{1}^{p}.u_{2}^{q}\ v^{p,q}$ pour $v^{p,q}\in V^{p,q}$. $\phi $ est un isomorphisme de ${\bf C}[u_{1},u_{2}]$-modules et donc $\xi_{{\bf A}^2}(V,F^{\bullet}_{1},F^{\bullet }_{2})$ est isomorphe au faisceau cohérent obtenu à partir du ${\bf C}[u_{1},u_{2}]$-module $\oplus_{(p,q) \in {\bf Z}^{2}} V^{p,q}$ i.e ${\mathcal O}_{{\bf A}^2}\otimes_{\bf C}V={\mathcal O}_{{\bf A}^2}^{dim_{\bf C}V}$. L'isomorphisme de modules $\phi$ induit un isomorphisme de faisceaux cohérents que nous noterons encore $\phi $ :
\begin{center}
$\phi :{\mathcal O}_{{\bf A}^2}\otimes_{\bf C}V \rightarrow \xi_{{\bf A}^2}(V,F^{\bullet}_{1},F^{\bullet }_{2})$.
\end{center}

\subsection{Construction du fibré de Rees sur ${\bf P}^2$ à partir de trois filtrations}

Recouvrons ${\bf P}^{2}=\text{Proj} \ {\bf C}[u_{1},u_{2},u_{3}]$ par les trois cartes affines standard $U_{k}=\{(u_{1},u_{2},u_{3}) \in {\bf P}^{2}, u_{k}\neq 0\}={\bf A}^{2}_{ij}=\text{Spec} \ {\bf C}[\frac{u_{i}}{u_k},\frac{u_{j}}{u_k}]$ où $i,j,k$ sont deux à deux distincts, $\{i,j,k\}=\{1,2,3\}$ et $i<j$. Sur les trois cartes, on effectue la construction du fibré de Rees associé à deux filtrations. Pour les trois triplets $(i,j,k)$ considérés au-dessus on construit le fibré de Rees $\xi_{{\bf A}^{2}_{ij}}(V,F^{\bullet}_{i},F^{\bullet }_{j})$ sur $U_{k}={\bf A}^{2}_{ij}$. D'après la section précédente, chacun de ces fibrés est isomorphe au fibré trivial sur ${\mathcal O}_{U_k}\otimes V$ sur $U_k$ par les isomorphismes $\phi_{k}:{\mathcal O}_{U_k} \rightarrow \xi_{{\bf A}^{2}_{ij}}(V,F^{\bullet}_{i},F^{\bullet }_{j})$. On peut recoller ces fibrés triviaux sur les cartes affines pour former un fibré sur ${{\bf P}^2}$ comme le montre la proposition suivante.

\begin{proposition-definition}\label{recolledef}
Le recollement des fibrés vectoriels triviaux ${\mathcal O}_{U_k}\otimes V$ où $k \in \{0,1,2\}$ par l'intermédiaire des restrictions des isomorphismes $\phi_k$ forment un fibré vectoriel sur ${\bf P}^2$. Le fibré ainsi obtenu est appelé fibré de Rees associé à l'espace vectoriel trifiltré  $(V,F_{0}^{\bullet },F_{1}^{\bullet },F_{2}^{\bullet })$ et est noté ${\xi}_{{\bf P}^2}(V,F_{0}^{\bullet },F_{1}^{\bullet },F_{2}^{\bullet })$.

\end{proposition-definition}

\begin{preuve} On utilise les restrictions des isomorphismes trivialisants des fibrés de Rees sur les plans affines décrits plus haut pour former des cocycles définissant un fibré. Introduisons quelques notations. Les intersections deux à deux des ouverts affines standards $U_{k}$ où $k \in \{0,1,2 \}$ seront notées $U_{kl}=U_{k}\cap U_{l}$ où $k,l \in \{0,1,2 \}$ et $k < l$. $j_{kl}$ sera l'inclusion $j_{kl}: U_{kl} \hookrightarrow U_{k}$ et $j_{lk}$ sera l'inclusion $j_{lk}: U_{kl} \hookrightarrow U_{l}$. De m\^eme, on a les trois inclusions $j_{ijk} : U_{ijk} \hookrightarrow U_{ij}$. 
Pour $l \neq k$, on peut supposer par exemple que $l=j$ où $(i,j,k)$ est un triplet de la forme décrite plus haut. On veut construire un morphisme de transition à partir des isomorphismes $\phi_{k}:{\mathcal O}_{U_k} \rightarrow \xi_{{\bf A}^{2}_{ij}}(V,F^{\bullet}_{i},F^{\bullet }_{j})$ et $\phi_{j}:{\mathcal O}_{U_j} \rightarrow \xi_{{\bf A}^{2}_{ik}}(V,F^{\bullet}_{i},F^{\bullet }_{k})$ sur $U_{jk}$. Or $ \xi_{{\bf A}^{2}_{ij}}(V,F^{\bullet}_{i},F^{\bullet }_{j})\vert_{U_{jk}}= \xi_{{\bf A}^{2}_{ik}}(V,F^{\bullet}_{i},F^{\bullet }_{k})\vert_{U_{jk}}$ comme sous faisceaux de $j^{*}_{jk}({\mathcal O}_{U_{ijk}}\otimes V)$. En effet ces deux faisceaux sont égaux au faisceau $\xi_{{\bf A}^1}(V,F_{i}^{\bullet })\otimes {\mathcal O}_{{\bf G}_{m}}$ sur $U_{jk}=Spec \ {\bf C}[u,v,v^{-1}]$, où ${\bf G}_{m}=\text{Spec} \ {\bf C}[v,v^{-1}]$. On peut ainsi définir $\phi_{kl}=\phi_{l}^{-1}\circ \phi_{k}\vert_{U_{kl}}$ comme fonction de transition entre les fibrés triviaux sur les cartes $U_{k} $ et $U_{l}$. La condition de cocycle est automatiquement vérifiée car $\phi_{k}\vert_{U_{ijk}}(\xi_{{\bf A}^{2}_{ij}}(V,F^{\bullet}_{i},F^{\bullet }_{j}))={\mathcal O}_{U_{ijk}}\otimes V$. Ainsi sur $U_{ijk}$ : $\phi_{ij} \circ \phi_{jk} \circ \phi_{kl}=id$.
Reste à montrer que la construction ne dépend pas des bigraduations $V=\oplus_{p,q}V^{p,q}$ choisie pour chaque construction $(\xi_{{\bf A}^{2}_{ij}}(V,F^{\bullet}_{i},F^{\bullet }_{j})$ sur les cartes affines $U_{ij}$. Supposons que l'on ait une autre graduation $V=\oplus_{p,q}U^{p,q}$, alors prendre le bigradué associé aux filtrations $F^{\bullet }_{i}$ et $F^{\bullet}_{j}$, donne des isomorphismes $V^{p,q} \cong U^{p,q}$. Ces isomorphismes induisent des isomorphismes au niveau des modules de Rees puis des fibrés de Rees sur les cartes affines $U_{k}$ pour tous les triplets $i,j,k$. Ainsi, changer les bigraduations associées aux paires de filtrations revient à changer de trivialisations locales du fibré localement libre sur ${\bf P}^2$ et conduit donc à un fibré isomorphe. Le fibré de Rees ${\xi}_{{\bf P}^2}(V,F_{0}^{\bullet },F_{1}^{\bullet },F_{2}^{\bullet })$ est bien déterminé à isomorphisme prés à partir de $(V,F_{0}^{\bullet },F_{1}^{\bullet },F_{2}^{\bullet })$.     
\end{preuve}







\textit{Exemple :}
  Donnons la description explicite du fibré de Rees associé à un espace vectoriel de dimension $1$ muni de trois filtrations. Cet exemple est important car on essaiera toujours par la suite de se ramener à une somme directe d'espaces vectoriels de dimension 1 munis de trois filtrations, la construction précédente se comportant bien pour les sommes directes. Soit $V={\bf C}$ muni de trois filtrations. On a forcément $(V,F_{0}^{\bullet },F_{1}^{\bullet },F_{2}^{\bullet })=(V ,Dec^{r}Triv^{\bullet},Dec^{p}Triv^{\bullet},Dec^{q}Triv^{\bullet})$ où $r,p,q$ sont les rangs respectifs auxquels $F^{\bullet}_{0}$, $F^{\bullet}_{1}$ et $F^{\bullet}_{2}$ sautent. ${\xi}_{{\bf P}^2}(V,Dec^{p}Triv^{\bullet},Dec^{q}Triv^{\bullet},Dec^{r}Triv^{\bullet})$ est un fibré en droites sur ${\bf P}^2$ donc de la forme $ O_{{\bf P}^2}(D)$. Explicitons le diviseur $D$. Pour cela notons $D_{k}$ le complémentaire de $U_{k}$ dans ${\bf P}^2$ (la droite a l'infini associée à l'ouvert affine $U_k$), c'est à dire : pour $k \in \{0,1,2\}$ : $U_{k}=\{ (u_{0}:u_{1}:u_{2}) \in {{\bf P}^2} \vert u_{k}\neq 0 \} \cong {\bf A}^{2}_{ij}$ et $D_{k}=\{ (u_{0}:u_{1}:u_{2}) \in {{\bf P}^2} \vert u_{k}= 0 \} \cong {\bf P}^{1}$. La construction dans le cas d'un fibré de rang $1$ est simplifiée par le fait que l'on peut évidemment trouver des trivialisations sur les ouverts compatibles à toutes les filtrations. Une section méromorphe $s$ du fibré est donnée par (on note $v$ un vecteur engendrant la droite $V$) : $ (\frac{u_0}{u_2})^{-r}.(\frac{u_1}{u_2})^{-p}v $ sur $U_{2}$, $ (\frac{u_1}{u_0})^{-p}.(\frac{u_2}{u_0})^{-q}v $ sur $U_{0}$ et $ (\frac{u_0}{u_1})^{-r}.(\frac{u_2}{u_1})^{-q}v $ sur $U_{1}$. D'où : $D=(s)=-rD_{0}-pD_{1}-qD_{2}$.

\begin{lemma}\label{dirsumxi}Soit $(V_{i},{F_{0}^{\bullet }}_{i},{F_{1}^{\bullet }}_{i},{F_{2}^{\bullet }}_{i})$ une famille finie d'espaces vectoriels munis de trois filtrations, alors :\begin{center}
 ${\xi}_{{\bf P}^2}(\oplus_{i}(V_{i},{F_{0}^{\bullet }}_{i},{F_{1}^{\bullet }}_{i},{F_{2}^{\bullet }}_{i})) \cong \oplus_{i}{\xi}_{{\bf P}^2}(V_{i},{F_{0}^{\bullet }}_{i},{F_{1}^{\bullet }}_{i},{F_{2}^{\bullet }}_{i}) $. 
  ${\xi}_{{\bf P}^2}(\otimes_{i}(V_{i},{F_{0}^{\bullet }}_{i},{F_{1}^{\bullet }}_{i},{F_{2}^{\bullet }}_{i})) \cong \otimes_{i}{\xi}_{{\bf P}^2}(V_{i},{F_{0}^{\bullet }}_{i},{F_{1}^{\bullet }}_{i},{F_{2}^{\bullet }}_{i}) $.
\end{center}
\end{lemma}
 
\begin{preuve}
Rappelons (\cite{har}) que si $M$ et $N$ sont deux $A$-modules, et que l'on note par une tilde la construction du faisceau cohérent sur $Spec \ A$ associée à un $A$-module, alors : $ (M \oplus N)^{\sim} \cong M^{\sim } \oplus N^{\sim } $ et $(M \otimes_{A} N)^{\sim}  \cong M^{\sim} \otimes_{Spec \ A} N^{\sim}$. Les constructions des fibrés de Rees sur les plans affines sur les ouverts standards associées à la somme directe et au produit tensoriel se déduisent donc des constructions sur chacun des facteurs. Les trivialisations sur ces ouverts sont données par chacunes des trivialisations correspondantes aux espaces vectoriels trifiltrés. Les fonctions de transitions du fibré associé à la somme directe et du fibré associé au produit tensoriel sont donc données par blocs correspondants aux somme directes et par produit tensoriel des fonctions de transition. Ainsi le fibré de Rees de la somme directe est bien la somme de Whitney des fibrés de Rees associés aux espaces vectoriels en somme directe et le fibré de Rees associé au produit tensoriel est le produit tensoriel des fibrés de Rees de chacun des facteurs.
\end{preuve}

\subsection{Etude du fibré de Rees}

L'étude du fibré de Rees associé à un espace vectoriel muni de trois filtrations (décroissantes et complètes) et de ses invariants va nous permettre de déceler à quel point les filtrations sont loins d'\^etre dans la position la plus simple, celle où elles sont simultanéement scindées. Rappelons que les trois filtrations $(F_{0}^{\bullet },F_{1}^{\bullet },F_{2}^{\bullet })$ sont simultanéement scindées s'il existe des sous-espaces vectoriels $V^{p,q}$ tels que :
\begin{center}$V =\oplus_{p,q}V^{p,q}$ (où la somme est finie) et 
$\left\lbrace \begin{array}{l}
         F_{0}^{r}=\oplus_{p'+q' \leq r}\,V^{p',q'},\\
         F_{1}^{p}=\oplus_{p'\leq p,q'}\,V^{p',q'},\\
         F_{2}^{q}=\oplus_{p',q' \leq q}\,V^{p',q'}.
\end{array}
\right.$.
\end{center}

\begin{lemma}\label{scin}
Soit $V$ un espace vectoriel muni de trois filtrations $(F_{0}^{\bullet },F_{1}^{\bullet },F_{2}^{\bullet })$ simultanéement scindées et de scindement $V=\oplus_{p,q}V^{p,q}$. Alors :
\begin{center}
${\xi}_{{\bf P}^2}(V,F_{0}^{\bullet },F_{1}^{\bullet },F_{2}^{\bullet })=\oplus_{p,q}{\xi}_{{\bf P}^2}(V^{p,q},Dec^{r}Triv^{\bullet },Dec^{p}Triv^{\bullet },Dec^{q}Triv^{\bullet })$.
\end{center}
\end{lemma}

\begin{preuve}
Il suffit d'écrire que $(V,F_{0}^{\bullet },F_{1}^{\bullet },F_{2}^{\bullet })=\oplus_{p,q}(V^{p,q},{F_{0}^{\bullet }}_{ind},{F_{1}^{\bullet }}_{ind},{F_{2}^{\bullet }}_{ind})$ et que sur chaque "bloc" $V^{p,q}$, ${F_{i}^{\bullet }}_{ind}$ est une filtration de niveau $1$, c'est à dire qu'il existe $k_{i}$ tel que $F^{k_i}_{i\  ind}=V^{p,q}$ et $F^{k_i}_{i\  ind}=\{0\}$. Sur $V^{p,q}$ on a $k_{i}=p$. Et $F^{\bullet }_{i\, ind}$ est donc une filtration triviale décalée $Dec^{p}Triv^{\bullet }$. On applique ensuite le lemme \ref{dirsumxi}. 
\end{preuve}

Nous n'avons pas de moyen, pour étudier le fibré de Rees associé à trois filtrations ${\xi}_{{\bf P}^2}(V,F_{0}^{\bullet },F_{1}^{\bullet },F_{2}^{\bullet })$, d'écrire des suites exactes en le scindant par l'une des filtrations pour faire appara\^{\i}tre des fibrés de rang plus bas construits à partir de filtrations sont en positions plus simples (on pense ici à la filtration par le poids pour une structure de Hodge mixte et aux structures pures sur chacun des gradué par le poids). Ceci vient du fait que si l'on peut scinder deux filtrations en m\^eme temps, ce n'est pas en général possible pour trois filtrations. Pour pouvoir étudier ${\xi}_{{\bf P}^2}(V,F_{0}^{\bullet },F_{1}^{\bullet },F_{2}^{\bullet })$ on va lui associer un fibré que l'on pourra démonter suivant les sous-espaces associés aux différentes filtrations.

Eclatons ${\bf P}^2$ en $(0:0:1)$, on notera $e: \widetilde{{\bf P}^2} \rightarrow {\bf P}^2$ le morphisme correspondant à l'éclatement et $E=e^{-1}(0)$ le diviseur exeptionnel. L'éclatement est décrit dans l'ouvert $U_{0}={\bf A}^{2}_{12}=Spec\,{\bf C}[u,v]$ par les morphismes $e_{i}: Spec\,B_{i} \rightarrow Spec\,A$ où $i \in \{1,2\}$ correspondants aux morphismes d'anneaux suivants (on note les morphismes d'anneaux de la m\^eme façon) : $e_{0}: A={\bf C}[u,v] \rightarrow B_{0}={\bf C}[x,y]$ tel que $e_{0}(u,v)=(uv,v)$ et  $e_{1}: A={\bf C}[u,v] \rightarrow B_{1}={\bf C}[z,t]$ tel que $e_{0}(u,v)=(u,uv)$. $Spec\,B_{i}$ est noté $U_{0}^{i}$. On fait la construction du fibré de Rees sur $Spec\,B_{i}$ associée au $B_{i}$-module $R^{2}(V,F_{i}^{\bullet },Triv^{\bullet })$.\\

\textit{Remarque :}
Sur l'ordre des filtrations. Le but du travail sur le fibré ${\xi}_{{\bf P}^2}(V,F_{0}^{\bullet },F_{1}^{\bullet },F_{2}^{\bullet })$ est de conna{\^\i}tre les positions relatives des filtrations dans le cas où l'espace vectoriel trifiltré est un structure de Hodge mixte. Dans ce cas, on prendra $F_{0}^{\bullet }=(W_{\bullet })^{\bullet }$ la filtration croissante associée à la filtration par le poids, et $F_{1}^{\bullet }=F_{}^{\bullet }$, $F_{2}^{\bullet }={\overline F}_{}^{\bullet }$ la filtration de Hodge et sa conjuguée par rapport à la structure réèlle sous-jacente. On va s'intéresser aux quotients par la filtration par le poids qui sont des structures de Hodge pures. C'est pourquoi dans la construction du fibré associé à ${\xi}_{{\bf P}^2}(V,F_{0}^{\bullet },F_{1}^{\bullet },F_{2}^{\bullet })$, on distingue $F_{0}^{\bullet }$ pour pouvoir quotienter par des sous-espaces vectoriels associés à cette filtration alors que les trois filtrations jouent des r\^oles symmétriques dans la construction du fibré de Rees.

\begin{proposition-definition}
On peut utiliser les trivialisations locales des fibrés de Rees associés aux paires de filtrations sur les quatres plans affines qui recouvrent $\widetilde{{\bf P}^2}$ pour obtenir un fibré vectoriel unique à isomorphisme prés. Il sera appelé le fibré vectoriel sur $\widetilde{{\bf P}^2}$ associé à  ${\xi}_{{\bf P}^2}(V,F_{0}^{\bullet },F_{1}^{\bullet },F_{2}^{\bullet })$, il sera noté ${\xi}_{\widetilde{{\bf P}^2}}(V,F_{0}^{\bullet },F_{1}^{\bullet },F_{2}^{\bullet },Triv^{\bullet })$.
\end{proposition-definition}

\begin{preuve}
La preuve est la m\^eme que pour la construction de ${\xi}_{{\bf P}^2}(V,F_{0}^{\bullet },F_{1}^{\bullet },F_{2}^{\bullet })$, on trivialise localement sur les quatres ouverts par l'intermédiaire des isomorphismes $\alpha_{i}$ exhibés plus haut ce qui permet d'avoir automatiquement les conditions de cocycles et d'utiliser de recoller.
\end{preuve}







\begin{lemma}\label{dirsumxitilde}Soit $(V_{i},{F_{0}^{\bullet }}_{i},{F_{1}^{\bullet }}_{i},{F_{2}^{\bullet }}_{i})$ une famille finie d'espaces vectoriels munis de trois filtrations, alors :\begin{center}
${\xi}_{\widetilde{{\bf P}^2}}(\oplus_{i}(V_{i},{F_{0}^{\bullet }}_{i},{F_{1}^{\bullet }}_{i},{F_{2}^{\bullet }}_{i},Triv^{\bullet })) \cong \oplus_{i}{\xi}_{\widetilde{{\bf P}^2}}(V_{i},{F_{0}^{\bullet }}_{i},{F_{1}^{\bullet }}_{i},{F_{2}^{\bullet }}_{i},Triv^{\bullet})$. \\
${\xi}_{\widetilde{{\bf P}^2}}(\otimes_{i}(V_{i},{F_{0}^{\bullet }}_{i},{F_{1}^{\bullet }}_{i},{F_{2}^{\bullet }}_{i},Triv^{\bullet })) \cong \otimes_{i}{\xi}_{\widetilde{{\bf P}^2}}(V_{i},{F_{0}^{\bullet }}_{i},{F_{1}^{\bullet }}_{i},{F_{2}^{\bullet }}_{i},Triv^{\bullet})$.  
\end{center}
\end{lemma}
\begin{preuve}
La preuve est similaire à celle du lemme \ref{dirsumxi}.   

\end{preuve}

\begin{definition}
Une filtration décroissante et complète $F^{\bullet }$ est dite positive si le plus petit entier  $p$ tel que $F^{k}=V$ pour tout $k \leq p$ est positif.
\end{definition}

Les fibrés ainsi définis sont liés par le lemme suivant :

\begin{lemma}Soit $(V,F^{\bullet }_{0},F^{\bullet }_{1},F^{\bullet }_{2})$ un espace vectoriel muni de trois filtrations. Si $F_{1}^{\bullet}$ et $F_{2}^{\bullet}$ sont positives on a un morphisme injectif de ${\xi}_{\widetilde{{\bf P}^2}}(V,F_{0}^{\bullet },F_{1}^{\bullet },F_{2}^{\bullet },Triv^{\bullet })$ vers $e^{*}{\xi}_{{\bf P}^2}(V,F_{0}^{\bullet },F_{1}^{\bullet },F_{2}^{\bullet })$.
\end{lemma}

\begin{preuve}
Démontrons ceci sur chacun des ouverts de $\widetilde{{\bf P}^2}$. Notons  $U_{1}'\cong U_{1}$ et $U_{2}'\cong U_{2}$ les images inverses de $U_{1}$ et $U_{2}$ par $e$. Sur les ouverts $U_{1}'$ et $U_{2}'$ il n'y a rien a montrer. Tout se passe au dessus de la carte $U_{0}$ dans ${\bf P}^2$ (rappelons qu'un recouvrement affine de $e^{-1}(U_{0}$ est donné par $U_{0}^{0}$ et $U_{0}^{1}$). Comme on l'a vu plus haut, pour un $A$-module $M$ et un morphisme d'anneaux $A \rightarrow B$, on a $f^{*}({\tilde M}) \cong \widetilde{M \otimes_{A} B }$. Par exactitude de $^\sim$ il suffit donc de montrer qu'on a une flèche de $B_{i}$-modules $RR(V,F_{i}^{\bullet },Triv^{\bullet }) \rightarrow  RR(V,F_{i}^{\bullet },F_{2}^{\bullet }) \otimes_{A} B_{i}$ pour $i \in \{1,2\}$ injective. Il suffit de le montrer pour $i=1$.\\
 $RR(V,F_{1}^{\bullet },F_{2}^{\bullet }) \otimes_{A} B_{1}$ est engendré par les éléments de la forme $u^{-p}.v^{-q}.a_{p,q} \otimes_{A}P(x,y)$ où $P \in B_{1}$ et $a_{p,q} \in F_{1}^{p } \cap F_{2}^{q }$, or $u^{-p}.v^{-q}.a_{p,q} \otimes_{A}P(x,y)=a_{p,q} \otimes_{A} e_{1}(u^{-p}).e_{1}(v^{-q}).P(x,y)=a_{p,q} \otimes_{A} (x.y)^{-p}.y^{-q}.P(x,y)=a_{p,q} \otimes_{A} x^{-p}.y^{-p-q}.P(x,y)$. Prenons comme flèche le morphisme injectif : $R^{2}(V,F_{1}^{\bullet },Triv^{\bullet }) \rightarrow  R^{2}(V,F_{1}^{\bullet },F_{2}^{\bullet }) \otimes_{A} B_{1}$ donnée par $x^{-p}.a_{p} \mapsto a_{p} \otimes_{A}x^{-p}$ où $a_{p} \in  F_{1}^{p }$.
\end{preuve}

\textit{Remarque :}
Lorsque le fibré sera associé à une struture de Hodge $F_{1}^{\bullet }=F_{}^{\bullet }$ et $F_{2}^{\bullet }={\overline F}_{}^{\bullet }$ et les filtrations au rang $p$ sont données par des $p$-formes et sont donc positives.   
Notons $\mathcal F $ le faisceau quotient de $e^{*}{\xi}_{{\bf P}^2}(V,F_{0}^{\bullet },F_{1}^{\bullet },F_{2}^{\bullet })$ par ${\xi}_{\widetilde{{\bf P}^2}}(V,F_{0}^{\bullet },F_{1}^{\bullet },F_{2}^{\bullet },Triv^{\bullet })$, on a la suite exacte suivante : 
\begin{center}
$0 \rightarrow {\xi}_{\widetilde{{\bf P}^2}}(V,F_{0}^{\bullet },F_{1}^{\bullet },F_{2}^{\bullet },Triv^{\bullet }) \rightarrow e^{*}{\xi}_{{\bf P}^2}(V,F_{0}^{\bullet },F_{1}^{\bullet },F_{2}^{\bullet }) \rightarrow {\mathcal F} \rightarrow 0$.
\end{center}
La construction de ${\xi}_{\widetilde{{\bf P}^2}}(V,F_{0}^{\bullet },F_{1}^{\bullet },F_{2}^{\bullet },Triv^{\bullet })$ est légitimée par le lemme suivant. Il permet de construire des suites exacte de fibrés sur $\widetilde{{\bf P}^2}$ associés aux sous-espaces et espaces quotients de la filtration $F^{\bullet }_{0}$ et donc de réduire l'étude de ${\xi}_{\widetilde{{\bf P}^2}}(V,F_{0}^{\bullet },F_{1}^{\bullet },F_{2}^{\bullet },Triv^{\bullet })$ puis de ${\xi}_{{\bf P}^2}(V,F_{0}^{\bullet },F_{1}^{\bullet },F_{2}^{\bullet })$ à celle des fibrés de Rees qui sont des fibrés en droite et donc de calculer les invariants topologiques du fibré de Rees associé à trois filtrations.

\begin{lemma}\label{estilde}
Soit $(V,F_{0}^{\bullet },F_{1}^{\bullet },F_{2}^{\bullet },Triv^{\bullet })$ un espace vectoriel muni de trois filtrations complètes et décroissantes, si $V'$ est un sous-espace vectoriel de $V$ donné par un des termes de la filtration $F_{0}^{\bullet }$(i.e il existe $p$ tel que $V'=F_{0}^{p}$), alors, en notant par des $'$ les sous-objets et les filtrations
induites sur $V'$ et par des $''$ les objets quotients et les filtrations quotient sur $V''=V/V'$, on a la suite exacte :
\begin{center}
\scalebox{0.9}[1]{$0 \rightarrow {\xi}_{\widetilde{{\bf P}^2}}(V',{F'}_{0}^{\bullet },{F'}_{1}^{\bullet },{F'}_{2}^{\bullet },Triv^{\bullet }) \rightarrow {\xi}_{\widetilde{{\bf P}^2}}(V,F_{0}^{\bullet },F_{1}^{\bullet },F_{2}^{\bullet },Triv^{\bullet }) \rightarrow {\xi}_{\widetilde{{\bf P}^2}}(V'',{F''}_{0}^{\bullet },{F''}_{1}^{\bullet },{F''}_{2}^{\bullet },Triv^{\bullet }) \rightarrow 0$.} 
\end{center}
\end{lemma}

\begin{preuve}
La suite est évidement exacte sur chacun des quatres ouverts affines decrits plus haut. Le recollement est rendu possible par le fait que sur chaque intersection triple deux filtrations différentes sont en jeu et donc toutes les trivialisations sont compatibles.
\end{preuve}

\textit{Remarque :}
L'idée donnée dans la démonstration du lemme nous donne m\^eme un résultat plus fort : tous les quotients de $(V,F_{0}^{\bullet },F_{1}^{\bullet },F_{2}^{\bullet },Triv^{\bullet })$ gradués par $F^{\bullet }_{0}$ forment des sous-fibrés de ${\xi}_{\widetilde{{\bf P}^2}}(V,F_{0}^{\bullet },F_{1}^{\bullet },F_{2}^{\bullet },Triv^{\bullet })$ et leurs quotients sont de m\^eme des sous-fibrés. 

Le lemme suivant est l'analogue du lemme \ref{scin} :

\begin{lemma}
Soit $V$ un espace vectoriel muni de trois filtrations $(F_{0}^{\bullet },F_{1}^{\bullet },F_{2}^{\bullet })$ scindées. Alors :
\begin{center}
${\xi}_{\widetilde{{\bf P}^2}}(V,F_{0}^{\bullet },F_{1}^{\bullet },F_{2}^{\bullet },Triv^{\bullet})=\oplus_{p,q}{\xi}_{\widetilde{{\bf P}^2}}(V^{p,q},Dec^{r}Triv^{\bullet },Dec^{p}Triv^{\bullet },Dec^{q}Triv^{\bullet },Triv^{\bullet})$.
\end{center}
\end{lemma}

\begin{preuve}
La preuve se calque sur celle du lemme \ref{scin}.
\end{preuve}

\subsection{Calcul explicites des invariants des fibrés de Rees}

Le but de cette section est de calculer le caractère de Chern du fibré ${\xi}_{{\bf P}^2}(V,F_{0}^{\bullet },F_{1}^{\bullet },F_{2}^{\bullet })$ associé à un
espace vectoriel trifiltré $(V,F_{0}^{\bullet },F_{1}^{\bullet },F_{2}^{\bullet })$ afin de voir en quoi ce fibré différe du fibré trivial obtenu pour la construction associée avec trois filtrations simultanéement scindées provenant des gradués donnés par les trois filtrations. Pour cela, on utilise la décomposition de la section précédente pour calculer les caractères de Chern $\text{ch}({\mathcal F})$ et $\text{ch}({\xi}_{\widetilde{{\bf P}^2}}(V,F_{0}^{\bullet },F_{1}^{\bullet },F_{2}^{\bullet },Triv^{\bullet }))$.

\begin{prop}
 Le calcul du caractère de Chern $\emph{ch}({\xi}_{\widetilde{{\bf P}^2}}(V,F_{0}^{\bullet },F_{1}^{\bullet },F_{2}^{\bullet },Triv^{\bullet }))$ peut toujours se ramener au calcul des caractères de Chern $\emph{ch}({\xi}_{\widetilde{{\bf P}^2}}({\bf C},Dec^{r}Triv^{\bullet },Dec^{p}Triv^{\bullet },Dec^{q}Triv^{\bullet },Triv^{\bullet }))$ pour $(r,p,q) \in {\bf Z}^3$. De plus, en notant  $\eta_{D} \in H^{2}({\widetilde{{\bf P}^2}},{\bf Z})={\bf Z}$ le dual de Poincaré du diviseur $D$ et ${\tilde w}^{4}$ une forme qui engendre $H^{4}({\widetilde{{\bf P}^2}},{\bf Z})={\bf Z}$ :\\
$\emph{ch}({\xi}_{\widetilde{{\bf P}^2}}({\bf C},Dec^{r}Triv^{\bullet },Dec^{p}Triv^{\bullet },Dec^{q}Triv^{\bullet },Triv^{\bullet }))=1+r {\eta}_{{\tilde D}_{0}}+ p{\eta}_{{\tilde D}_{1}}+q{\eta}_{{\tilde D}_{2}}+{\frac{1}{2}}(r^2+2rp+2rq){\tilde w}^{4}$
On a ainsi la relation :\\
$\emph{ch}({\xi}_{\widetilde{{\bf P}^2}}(V,F_{0}^{\bullet },F_{1}^{\bullet },F_{2}^{\bullet },Triv^{\bullet }))={\emph{dim}}_{\bf C}V+ \sum_{p,q,r}{\delta }_{p,q,r}.(r {\eta}_{{\tilde D}_{0}}+ p{\eta}_{{\tilde D}_{1}}+q{\eta}_{{\tilde D}_{2}}+{\frac{1}{2}}(r^2+2rp+2rq){\tilde w}^{4})$
où ${\delta }_{p,q,r}={\emph{dim}}_{\bf C}Gr_{F_{2}^{\bullet }}^{q}Gr_{F_{1}^{\bullet }}^{p}Gr_{F_{0}^{\bullet }}^{r}V$.
\end{prop}

Avant de démontrer la proposition, établissons le lemme suivant :

\begin{lemma}$ $\\
${\xi}_{{\bf P}^2}({\bf C},Dec^{r}Triv^{\bullet },Dec^{p}Triv^{\bullet },Dec^{q}Triv^{\bullet })\cong O_{{\bf P}^2}(rD_{0}+pD_{1}+qD_{2})$.\\ 
$c_{1}({\xi}_{{\bf P}^2}({\bf C},Dec^{r}Triv^{\bullet },Dec^{p}Triv^{\bullet },Dec^{q}Triv^{\bullet }))=r \eta_{D_{0}}+p\eta_{D_{1}}+q\eta_{D_{2}}$.\\${\xi}_{\widetilde{{\bf P}^2}}({\bf C},Dec^{r}Triv^{\bullet },Dec^{p}Triv^{\bullet },Dec^{q}Triv^{\bullet },Triv^{\bullet })\cong O_{\widetilde {{\bf P}^2}}(r{\tilde D}_{0}+p{\tilde D}_{1}+q{\tilde D}_{2})$.\\
$c_{1}({\xi}_{\widetilde{{\bf P}^2}}({\bf C},Dec^{r}Triv^{\bullet },Dec^{p}Triv^{\bullet },Dec^{q}Triv^{\bullet },Triv^{\bullet }))=r {\eta}_{{\tilde D}_{0}}+ p{\eta}_{{\tilde D}_{1}}+q{\eta}_{{\tilde D}_{2}}.$

\end{lemma}

\begin{preuve} La première assertion a été démontrée dans l'exemple qui suit la proposition-definition \ref{recolledef}. Pour distinguer plus bas leurs transformées strictes dans l'éclaté $\widetilde{{\bf P}^2}$ de ${\bf P}^2$, on garde les notations des $D_{i}$ bien que tous soient homologues dans ${\bf P}^2$ et  définissent des fibrés isomorphes. D'après \cite{grihar} par exemple, sur une variété compacte complexe $M$, la classe de Chern d'un fibré en droite de la forme $O_{M}(D)$ pour $D \in \text{Div}(M)$ est donnée par $c_{1}(O_{M}(D))= \eta_{D}$ où $\eta_{D} \in H^{2}_{DR}(M) $ est le dual de Poincaré de $D$. Ainsi :\\
$c_{1}({\xi}_{{\bf P}^2}(V,Dec^{r}Triv^{\bullet },Dec^{p}Triv^{\bullet },Dec^{q}Triv^{\bullet }))= \eta_{D_{0}}+\eta_{D_{1}}+\eta_{D_{2}}$
Notons ${\tilde D}_{i}$ la transformée stricte dans $\widetilde{{\bf P}^2}$ de $D_{i}$ donnée par l'application de l'éclatement $e$, et $E$ le diviseur exeptionnel. Le fibré en droites ${\xi}_{\widetilde{{\bf P}^2}}({\bf C},Dec^{r}Triv^{\bullet },Dec^{p}Triv^{\bullet },Dec^{q}Triv^{\bullet },Triv^{\bullet }))$ sur ${\widetilde{{\bf P}^2}}$ est de la forme $O_{\widetilde {{\bf P}^2}}(D)$ pour $[D] \in \text{Pic}({\widetilde {{\bf P}^2}})={\bf Z}\oplus {\bf Z}$. Un calcul analogue à celui fait dans l'exemple précité montre l'égalité $D=r{\tilde D}_{0}+p{\tilde D}_{1}+q{\tilde D}_{2}$ et donc : $$c_{1}({\xi}_{\widetilde{{\bf P}^2}}({\bf C},Dec^{r}Triv^{\bullet },Dec^{p}Triv^{\bullet },Dec^{q}Triv^{\bullet },Triv^{\bullet }))=r {\eta}_{{\tilde D}_{0}}+ p{\eta}_{{\tilde D}_{1}}+q{\eta}_{{\tilde D}_{2}}.$$

\end{preuve}
Démontrons la proposition :
\begin{preuve} Démontrons la première assertion. On notera toujours par $Triv^{\bullet }$ la filtration triviale induite par la filtration triviale sur des sous-espaces ou des espaces quotients de $V$. Soit $V$ un espace vectoriel muni de trois filtrations opposées $(F_{0}^{\bullet },F_{1}^{\bullet },F_{2}^{\bullet })$. Rappelons que pour une suite exacte de faisceaux cohérents de la forme
 $0 \rightarrow {\mathcal G}' \rightarrow {\mathcal G} \rightarrow {\mathcal G}'' \rightarrow 0$ la relation suivante est vérifiée : $ \text{ch}( {\mathcal G})= \text{ch}({\mathcal G}')+ \text{ch}({\mathcal G}'')$. Ainsi en filtrant $V$ par $F_{0}^{\bullet }$ et d'après la suite exacte exhibée dans le lemme \ref{estilde} :

$$\text{ch}({\xi}_{\widetilde{{\bf P}^2}}(V,F_{0}^{\bullet },F_{1}^{\bullet },F_{2}^{\bullet },Triv^{\bullet }))=\sum_{r}\text{ch}({\xi}_{\widetilde{{\bf P}^2}}(Gr_{F_{0}^{\bullet }}^{r}V,F_{0,ind_{n}}^{\bullet },F_{1,ind_{r}}^{\bullet },{F_{2,ind_{}r}^{\bullet }},Triv^{\bullet })).$$

Comme $F_{0,ind_{r}}^{\bullet }$ est de longueur $1$ sur $Gr_{F_{0}^{\bullet }}^{r}V$ c'est la filtration décalée de $r$ par rapport à la filtration triviale, donc :

\hspace{1.5cm}$ \text{ch}({\xi}_{\widetilde{{\bf P}^2}}(Gr_{F_{0}^{\bullet }}^{r}V,F_{0,ind_{r}}^{\bullet },F_{1,ind_{r}}^{\bullet },{F_{2,ind_{r}}^{\bullet }},Triv^{\bullet }))=$

\hspace{6cm}  $ \text{ch}({\xi}_{\widetilde{{\bf P}^2}}(Gr_{F_{0}^{\bullet }}^{r}V,Dec^{r}Triv^{\bullet},F_{1,ind_{r}}^{\bullet },{F_{2,ind_{}r}^{\bullet }},Triv^{\bullet })).$

Pour tout $r$, filtrons $Gr_{F_{0}^{\bullet }}^{r}$ par la filtration induite par $F_{1,ind_{r}}^{\bullet }$ sur cet espace. Il vient alors :
 
\hspace{1.5cm}$ \text{ch}({\xi}_{\widetilde{{\bf P}^2}}(Gr_{F_{0}^{\bullet }}^{r}V,Dec^{r}Triv^{\bullet},F_{1,ind_{r}}^{\bullet },{F_{2,ind_{}r}^{\bullet }},Triv^{\bullet }))=$ 

\hspace{3.5cm}$\sum_{p}  \text{ch}({\xi}_{\widetilde{{\bf P}^2}}(Gr_{F_{1,ind_{r}}^{\bullet }}^{p}Gr_{F_{0}^{\bullet }}^{r}V,Dec^{r}Triv^{\bullet},F_{1,ind_{r},ind_{p}}^{\bullet },{F_{2,ind_{r},ind_{p}}^{\bullet }},Triv^{\bullet })).$

Or $F_{1,ind_{r},ind_{p}}^{\bullet }$ est de longueur $1$ sur $Gr_{F_{1,ind_{r}}^{\bullet }}^{p}Gr_{F_{0}^{\bullet }}^{r}V=Gr_{F_{1}^{\bullet }}^{p}Gr_{F_{0}^{\bullet }}^{r}V$, ainsi :

\hspace{1.5cm}$ \text{ch}({\xi}_{\widetilde{{\bf P}^2}}(Gr_{F_{1,ind_{r}}^{\bullet }}^{p}Gr_{F_{0}^{\bullet }}^{r}V,Dec^{r}Triv^{\bullet},F_{1,ind_{r},ind_{p}}^{\bullet },{F_{2,ind_{r},ind_{p}}^{\bullet }},Triv^{\bullet }))=$

\hspace{5cm}$ \text{ch}({\xi}_{\widetilde{{\bf P}^2}}(Gr_{F_{1}^{\bullet }}^{p}Gr_{F_{0}^{\bullet }}^{r}V,Dec^{r}Triv^{\bullet},Dec^{p}Triv^{\bullet },{F_{2,ind_{r},ind_{p}}^{\bullet }},Triv^{\bullet })).$

Le m\^eme argument appliqué à la dernière filtration permet de montrer que :

$ \text{ch}({\xi}_{\widetilde{{\bf P}^2}}(V,F_{0}^{\bullet },F_{1}^{\bullet },F_{2}^{\bullet },Triv^{\bullet }))=$

\hspace{3cm} $\sum_{p,q,r}  \text{ch}({\xi}_{\widetilde{{\bf P}^2}}(Gr_{F_{2}^{\bullet }}^{q}Gr_{F_{1}^{\bullet }}^{p}Gr_{F_{0}^{\bullet }}^{r}V,Dec^{r}Triv^{\bullet },Dec^{p}Triv^{\bullet },Dec^{q}Triv^{\bullet },Triv^{\bullet })).$

Puis finalement que :\\

$ \text{ch}({\xi}_{\widetilde{{\bf P}^2}}(V,F_{0}^{\bullet },F_{1}^{\bullet },F_{2}^{\bullet },Triv^{\bullet }))=$

\hspace{0.5cm}$ \sum_{p,q,r} \,  \text{dim}_{\bf C}( Gr_{F_{2}^{\bullet }}^{q}Gr_{F_{1}^{\bullet }}^{p}Gr_{F_{0}^{\bullet }}^{r}V).\, \text{ch}({\xi}_{\widetilde{{\bf P}^2}}({\bf C},Dec^{r}Triv^{\bullet },Dec^{p}Triv^{\bullet },Dec^{q}Triv^{\bullet },Triv^{\bullet })).$\\

Ce qui démontre la première assertion.\\
La deuxième découle du lemme précédent. Pour un fibré en droite $\mathcal L$, on a la relation $ \text{ch}({\mathcal L})=\exp(c_{1}({\mathcal L}))$ c'est à dire sur une surface : $ \text{ch}({\mathcal L})=1+c_{1}({\mathcal L})+{\frac{1}{2}}c_{1}({\mathcal L})^{2}$. Les produits d'intersection dans ${\bf {P^2}}$ sont donnés par ${{\tilde D}_{0}}^{2}=1,\,{{\tilde D}_{1}}^{2}=0,\,{{\tilde D}_{2}}^{2}=0,\,E^{2}=-1$ et ${{\tilde D}_{0}}.{{\tilde D}_{1}}=1,\,{{\tilde D}_{0}}.{{\tilde D}_{2}}=1,\,E.{{\tilde D}_{1}}=1,\,E.{{\tilde D}_{2}}=1,\,{{\tilde D}_{1}}.{{\tilde D}_{2}}=0$. D'où :
$$ \text{ch}({\xi}_{\widetilde{{\bf P}^2}}({\bf C},Dec^{r}Triv^{\bullet },Dec^{p}Triv^{\bullet },Dec^{q}Triv^{\bullet },Triv^{\bullet }))=1+r {\eta}_{{\tilde D}_{0}}+ p{\eta}_{{\tilde D}_{1}}+q{\eta}_{{\tilde D}_{2}}+{\frac{1}{2}}(r^2+2rp+2rq){\tilde w}^{4}.$$


\end{preuve}
Reste à déterminer $ \text{ch}(\mathcal F)$. L'idée pour le faire est de voir que $\mathcal F $ ne dépend pas de la filtration $F_{0}^{\bullet }$. En effet le support de $\mathcal F $ est dans l'éclaté de la carte affine $U_{0}$ et $\calf $ peut \^etre déterminé par les restrictions des faisceaux  ${\xi}_{{\bf P}^2}(V,F_{0}^{\bullet },F_{1}^{\bullet },F_{2}^{\bullet })$ et ${\xi}_{{\bf P}^2}(V,F_{0}^{\bullet },F_{1}^{\bullet },F_{2}^{\bullet })$ à l'éclaté de cet ouvert affine. Or ces restrictions sont construites à partir uniquement des filtrations $F^{\bullet }_{1}$, $F^{\bullet }_{2}$ et $Triv^{\bullet }$ et ne dépendent pas de la filtration $F_{0}^{\bullet }$. Donc on peut écrire une suite exacte courte :\\
 $0 \rightarrow {\xi}_{\widetilde{{\bf P}^2}}(V,F^{\bullet },F_{1}^{\bullet },F_{2}^{\bullet },Triv^{\bullet }) \rightarrow e^{*}{\xi}_{{\bf P}^2}(V,F^{\bullet },F_{1}^{\bullet },F_{2}^{\bullet }) \rightarrow {\mathcal F} \rightarrow 0$ 
pour n'importe quelle filtration complète décroissante $F^{\bullet }$ de $V$ avec le m\^eme faisceau $\calf$. Ainsi, on va choisir une filtration $F^{\bullet }$ telle que $ \text{ch}({\xi}_{{\bf P}^2}(V,F^{\bullet },F_{1}^{\bullet },F_{2}^{\bullet }))$ soit facilemnt calculable. Pour cela il faut qu'il y ait un scindement commun de $V$ par les trois filtrations $F^{\bullet }$, $F_{1}^{\bullet }$ et $F_{2}^{\bullet }$ ce qui est toujours le cas lorsqu'il n'y a que deux filtrations, comme dans le cas où $F^{\bullet }=F_{1}^{\bullet }$ ou $F^{\bullet }=F_{2}^{\bullet }$ par exemple.

\begin{prop}
La caractéristique de Chern de $\calf$ est $ \emph{ch}(\calf)=\frac{1}{2}(p+q)^{2}{\tilde w}^{4}$
\end{prop}

\begin{preuve}
  On a le choix sur $F^{\bullet }$ pour effectuer le calcul. Celui-ci est facilité si l'on prend par exemple $F^{\bullet }=F_{1}^{\bullet }$. La suite exacte courte de faisceaux explicitée plus haut nous donne :

$ \text{ch}(\calf)= \text{ch}(e^{*}{\xi}_{{\bf P}^2}(V,F_{1}^{\bullet },F_{1}^{\bullet },F_{2}^{\bullet })) -\text{ch}({\xi}_{\widetilde{{\bf P}^2}}(V,F_{1}^{\bullet },F_{1}^{\bullet },F_{2}^{\bullet },Triv^{\bullet }))$.\\
D'après la proposition précédente, en posant ${\delta }_{p,q}= \text{dim}_{\bf C}Gr_{F_{2}^{\bullet }}^{q}Gr_{F_{1}^{\bullet }}^{p}V$ :

 $  \text{ch}({\xi}_{\widetilde{{\bf P}^2}}(V,F_{1}^{\bullet },F_{1}^{\bullet },F_{2}^{\bullet },Triv^{\bullet }))= \text{dim}_{\bf C}V+ \sum_{p,q}{\delta }_{p,q}.(p {\eta}_{{\tilde D}_{0}}+ p{\eta}_{{\tilde D}_{1}}+q{\eta}_{{\tilde D}_{2}}+{\frac{1}{2}}(p^2+2p^{2}+2pq)w)$.\\
D'un isomorphisme $V \cong Gr_{F_{2}^{\bullet }}^{q}Gr_{F_{1}^{\bullet }}^{p}V$ (rappelons qu'il en existe toujours un) et du lemme \ref{dirsumxi} on tire :
\begin{eqnarray*}
{\xi}_{{\bf P}^2}(V,F_{1}^{\bullet },F_{1}^{\bullet },F_{2}^{\bullet }) & \cong & {\xi}_{{\bf P}^2}(\oplus_{p,q}Gr_{F_{2}^{\bullet }}^{q}Gr_{F_{1}^{\bullet }}^{p}V,\oplus_{p,q}Dec^{p}Triv^{\bullet }_{p,q},\oplus_{p,q}{Dec^{p}Triv^{\bullet }}_{p,q},\oplus_{p,q}Dec^{q}Triv^{\bullet }_{p,q})\\
& \cong &  \oplus_{p,q}{\xi}_{{\bf P}^2}({\bf C},Dec^{p}Triv^{\bullet},Dec^{p}Triv^{\bullet},Dec^{q}Triv^{\bullet})^{dim_{\bf C}Gr_{F_{2}^{\bullet }}^{q}Gr_{F_{1}^{\bullet }}^{p}V} .
\end{eqnarray*}
Donc, si l'on note $w^{4}$ le générateur de $H^{4}({{\bf P}^2},{\bf Z})$ tel que $e^{*}(w^{4})={\tilde w}^{4}$ :\\

$ \text{ch}(e^{*}{\xi}_{{\bf P}^2}(V,F_{1}^{\bullet },F_{1}^{\bullet },F_{2}^{\bullet }))=e^{*} \text{ch}({\xi}_{{\bf P}^2}(V,F_{1}^{\bullet },F_{1}^{\bullet },F_{2}^{\bullet }))$
\begin{eqnarray*}
&=&e^{*} \text{ch}(\oplus_{p,q}{\xi}_{{\bf P}^2}({\bf C},Dec^{p}Triv^{\bullet},Dec^{p}Triv^{\bullet},Dec^{q}Triv^{\bullet})^{ \text{dim}_{\bf C}Gr_{F_{2}^{\bullet }}^{q}Gr_{F_{1}^{\bullet }}^{p}V})\\
&=&e^{*}( \text{dim}_{\bf C}V+\sum_{p,q}{\delta}_{p,q}(p {\eta}_{{ D}_{0}}+ p{\eta}_{{ D}_{1}}+q{\eta}_{{ D}_{2}}+{\frac{1}{2}}(p^2+p^{2}+q^{2}+2(p^{2}+pq+pq))w^{4})\\
&=& \text{dim}_{\bf C}V+\sum_{p,q}{\delta}_{p,q}(p {\eta}_{{\tilde D}_{0}}+ p{\eta}_{{\tilde D}_{1}}+q{\eta}_{{\tilde D}_{2}}+{\frac{1}{2}}(p^{2}+p^{2}+q^{2}+2(p^{2}+pq+pq)){\tilde w}^{4}).
\end{eqnarray*}
D'après la formule : $ \text{ch}({\xi}_{{\bf P}^2}({\bf C},Dec^{r}Triv^{\bullet},Dec^{p}Triv^{\bullet},Dec^{q}Triv^{\bullet})=1+r \eta_{D_{0}}+p\eta_{D_{1}}+q\eta_{D_{2}}$

\hspace{6.5cm} $+{\frac{1}{2}} (r^{2}+p^{2}+q^{2}+2(rp+rq+pq))w^{4}$).\\
\\
Ainsi :
$$ \text{ch}(\calf)=\frac{1}{2}\sum_{p,q}{\delta}_{p,q}(p+q)^{2}{\tilde w}^{4}.$$

\end{preuve}
Finalement, la caractéristique de Chern du fibré ${\xi}_{{\bf P}^2}(V,F_{0}^{\bullet },F_{1}^{\bullet },F_{2}^{\bullet })$ est donnée par :
\begin{prop}
$ \emph{ch}({\xi}_{{\bf P}^2}(V,F_{0}^{\bullet },F_{1}^{\bullet },F_{2}^{\bullet })=\emph{dim}_{\bf C}V$
\begin{eqnarray*}
\hspace{2cm}&+&\sum_{p,q,r}\Bigl[(\emph{dim}_{\bf C}Gr_{F_{2}^{\bullet }}^{q}Gr_{F_{1}^{\bullet }}^{p}Gr_{F_{0}^{\bullet }}^{r}V).\Bigl((r+p+q){w}^{2}+{\frac{1}{2}}(r^2+2rp+2rq){ w}^{4})\Bigr)\Bigr]\\
\hspace{2cm}&+&\sum_{p,q}\Bigr[(\emph{dim}_{\bf C}Gr_{F_{2}^{\bullet }}^{q}Gr_{F_{1}^{\bullet }}^{p}V).({\frac{1}{2}}(p+q)^{2}{ w}^{4})\Bigr]
\end{eqnarray*}
\end{prop}

\begin{preuve}
C'est une conséquence directe du calcul de $ \text{ch}(\calf)$ et de $ \text{ch}\,{\xi}_{\widetilde{{\bf P}^2}}(V,F^{\bullet }_{1},F_{1}^{\bullet },F_{2}^{\bullet },Triv^{\bullet })$.
\end{preuve}
Nous allons effectuer ce calcul dans le cas qui nous intéressera par la suite c'est à dire celui où les trois filtrations sont opposées.
\begin{definition}
Soit $V$ un espace vectoriel muni de trois filtrations $(F_{0}^{\bullet },F_{1}^{\bullet },F_{2}^{\bullet })$ complètes,finies et décroissantes. Ces filtrations sont dites opposées si $Gr_{F_{1}^{\bullet }}^{p}Gr_{F_{2}^{\bullet }}^{q}Gr_{F_{0}^{\bullet }}^{n}=0$ pour $p+q+n \neq 0$.
\end{definition}
\begin{corollaire}
\label{filtrscin}
Soient $(F_{0}^{\bullet },F_{1}^{\bullet },F_{2}^{\bullet })$ trois filtrations complètes, finies, décroissantes et opposées sur un espace vectoriel de dimension finie $V$, alors :

\hspace{1cm} $ \emph{ch}({\xi}_{{\bf P}^2}(V,F_{0}^{\bullet },F_{1}^{\bullet },F_{2}^{\bullet })= \emph{dim}_{\bf C}V$

\hspace{3cm} $+{\frac{1}{2}}\sum_{p,q}\Bigl[ \emph{dim}_{\bf C}Gr_{F_{2}^{\bullet }}^{q}Gr_{F_{1}^{\bullet }}^{p}V- \emph{dim}_{\bf C}Gr_{F_{2}^{\bullet }}^{q}Gr_{F_{1}^{\bullet }}^{p}Gr_{F_{0}^{\bullet }}^{-p-q}V \Bigr].(p+q)^{2}.w^{4}$.
\end{corollaire}

\begin{preuve}
On applique la formule de la proposition précédente. Tous les coefficients sont nuls si $p+q+n \neq 0$.  
\end{preuve}

\textit{Remarque :} Lorsque $(F_{0}^{\bullet },F_{1}^{\bullet },F_{2}^{\bullet })$ sont trois filtrations opposées sur $V$, on a :

\hspace{4cm} $\text{ch}_{1}({\xi}_{{\bf P}^2}(V,F_{0}^{\bullet },F_{1}^{\bullet },F_{2}^{\bullet }) )=0$.


\section{Application aux structures de Hodge mixtes : niveau de ${\bf R}$-scindement}

\subsection{Définition du niveau de ${\bf R}$-scindement d'une structure de Hodge mixte}

Une structure de Hodge mixte sur un espace $H_{\bf Q}$ consiste en les données suivantes :\\
                 \hspace*{.5cm}
{\bf (i)} \ Une filtration croissante $W_{\bullet }$ sur $H_{\bf Q}$ appelée filtration par le poids.\\
                \hspace*{.4cm}
{\bf (ii)} \ Une filtration décroissante appelée filtration de Hodge $F^{\bullet }$ sur $H=H_{\bf Q}\otimes_{\bf Q}{\bf C}$ satisfaisant la condition qui suit : la filtration $F^{\bullet }$ induit une structure de Hodge pure de poids $n$ sur
$Gr^{W}_{n}H=W_{n}\otimes_{\bf Q}{\bf C}/W_{n-1}\otimes_{\bf Q}{\bf C}$.

Ici la filtration $F^{\bullet }Gr^{W}_{n}H$ induite par $F^{\bullet }$ sur $Gr^{W}_{n}H$ est donnée par les quotients successifs pour tout entier $p$ : $F^{p }Gr^{W}_{n}H=(F^{p} \cap (W_{n}\otimes_{\bf Q}{\bf C})+W_{n-1}\otimes_{\bf Q}{\bf C}))/ W_{n-1}\otimes_{\bf Q}{\bf C}$.\\
Pour \^etre plus précis la condition {\bf (ii)} signifie que $F^{\bullet }Gr^{W}_{n}H$ et ${\overline F}^{\bullet }Gr^{W}_{n}H$ sont opposées sur $Gr^{W}_{n}H$ i.e $Gr^{p}_{F^{\bullet }Gr^{W}_{n}H}Gr^{q}_{{\overline F}^{\bullet }Gr^{W}_{n}H}H \neq 0$ si et seulement si $p+q \neq 0$ où ${\overline F}^{\bullet}$ est la filtration conjuguée à la filtration de Hodge par rapport à la structure réelle sous-jacente $H_{\bf R}=H_{\bf Q} \otimes_{\bf Q} {\bf R}$.

Nous n'utiliserons pas par la suite la structure rationnelle $H_{\bf Q}$. Il nous suffira de considérer la filtration par le poids sur $H_{\bf C}=H_{\bf Q} \otimes_{\bf Q} {\bf C}$ que nous noterons $W_{\bullet }$ alors que nous la notions $W_{\bf C} \otimes_{\bf Q} {\bf C}$. 
\textit{Remarque :}(Filtrations opposées).

Rappelons une propriété importante dans l'étude des structure de Hodge mixtes. Bien que les structures de Hodge mixte n'admettent pas en général de décomposition en une somme directe de "$(p,q)$-sous-espaces" comme les structures pures, c'est-à-dire de graduation compatible aux trois filtrations, les espaces canoniquement définis ci-dessous ont des propriétés de décomposition trés utiles :
$$ I^{p,q}=(F^{p} \cap W_{p+q}) \cap ({\overline F}^{q} \cap W_{p+q}+\sum_{i \geq 1} \, {\overline F}^{q-i}\cap W_{p+q-i-1}).$$   
    
\begin{lemma}\cite{del2} 

\item{(i)} $I^{p,q}={\overline I}^{q,p} \ \ mod \ W_{p+q-2}. $

\item{(ii)} $W_{m}=\oplus_{p+q \leq m}I^{p,q}$.

\item{(iii)} $F^{p}=\oplus_{i \geq p} \oplus_{q} I^{i,q}$.

\item{(iv)} La projection $W_{m} \rightarrow Gr_{m}^{W}H$ induit un isomorphisme pour $p+q=m$ de $I^{p,q}$ vers le sous-espace de Hodge $(Gr^{W}_{m}H)^{p,q}$. 
\end{lemma}

Quitte à changer les indices, cette décomposition canonique en sous-espaces $I^{p,q}$ donne deux bigraduations canoniques $H_{\bf C}=\oplus_{p,q}I^{p,q}$ associées aux paires de filtrations $(W_{\bullet },F^{\bullet })$ et $(W_{\bullet },{\overline F}^{\bullet })$. La première bigraduation, associée à $(W_{\bullet },F^{\bullet })$ est donnée directement par le lemme précédent, la deuxième associée à $(W_{\bullet },{\overline F}^{\bullet })$ est donnée par $W_{m}=\oplus_{p+q \leq m}{\overline I}^{p,q}$ et ${\overline F}^{q}=\oplus_{i \geq q} \oplus_{p} {\overline I}^{i,p}$. Ces deux bigraduations ne sont pas compatibles à moins que l'on ne soit dans le cadre de la définition suivante :
\begin{definition}
Une structure de Hodge mixte $H=(H_{\bf Q},W_{\bullet },F^{\bullet },{\overline F}^{\bullet })$ est dite {\bf R}-scindée si les espaces $I^{p,q}$ définis de façon canonique vérifient la condition : pour tous $p,q$    $I^{p,q}={\overline I}^{q,p}$.
\end{definition}
Dans le cas où la structure de Hodge mixte $H=(H_{\bf Q},W_{\bullet },F^{\bullet },{\overline F}^{\bullet })$ est scindée, on a : $F^{p}\cap {\overline F}^{q}=\oplus_{p'\geq p,q' \geq q}I^{p',q'}$. Ce qui signifie que les trois filtrations sont simultanéement scindées. Nous voulons étudier à quel point nous sommes éloignés en général de cette situation.

A une structure de Hodge mixte $H=(H_{\bf Q},W_{\bullet },F^{\bullet },{\overline F}^{\bullet })$ on peut associer les nombres entiers suivants :
$$\left\lbrace \begin{array}{l}
         h^{p,q}_{H}=dim_{\bf C}Gr_{{\overline F}^{\bullet }}^{q}Gr_{F^{\bullet }}^{p}Gr_{W_{\bullet }}^{p+q}H_{\bf C} \ \mathrm{pour \ tous}\ p,q \in {\bf Z}, \\
         t^{p,q}_{H}=dim_{\bf C}Gr_{{\overline F}^{\bullet }}^{q}Gr_{F^{\bullet }}^{p}H_{\bf C}\ \mathrm{ pour\ tous}\ p,q \in {\bf Z}.
         
\end{array}
\right.$$ 
Les entiers $h^{p,q}_{H}$ sont les nombres de Hodge classiques associés aux structures de Hodge mixtes. Dans le cas où la structure de Hodge mixte $H$ est ${\bf R}$-scindée, pour tous $p,q \in {\bf Z}$ on a l'égalité : $h^{p,q}_{H}=t^{p,q}_{H}$. Ce n'est pas vrai en général comme nous le verrons par la suite. A une structure de Hodge mixte $(H_{\bf Q},W_{\bullet },F^{\bullet},{\overline F}^{\bullet})$ on  associe le fibré de Rees sur ${\bf P}^2$ construit dans section précédente à partir du triplet de filtrations sur $H_{\bf C}$, $(W_{.}^{\bullet },F^{\bullet },{\overline F}^{\bullet })$. $W^{\bullet }_{.}$ est la filtration décroissante associée à la filtration croissante $W_{\bullet }$. Elle est définie par : pour tout $p \in {\bf Z}$,  $W_{.}^{p}=W_{-p}$. Ainsi à toute structure de Hodge mixte $(H_{\bf Q},W_{\bullet },F^{\bullet},{\overline F}^{\bullet})$ correspond un fibré vectoriel sur le plan projectif ${\xi}_{{\bf P}^2}(H_{\bf C},W_{.}^{\bullet },F^{\bullet },{\overline F}^{\bullet })$. Le calcul des invariants du fibré construit plus haut nous permettra de voir si les nombres de Hodge $t^{p,q}_{H}$ sont proches des nombres $h^{p,q}_{H}$, c'est à dire de voir si la structure de Hodge mixte est proche de la structure de Hodge mixte ${\bf R}$-scindée qui lui est associée \cite{catkapsch}.

\begin{prop}
Soit ${\xi}_{{\bf P}^2}(H_{\bf C},W_{.}^{\bullet },F^{\bullet },{\overline F}^{\bullet })$ le fibré sur ${\bf P}^2$ associé à la structure de Hodge mixte $(H_{\bf Q},W_{\bullet },F^{\bullet},{\overline F}^{\bullet})$, alors :
$$ \emph{ch}({\xi}_{{\bf P}^2}(H_{\bf C},W_{.}^{\bullet },F^{\bullet },{\overline F}^{\bullet }))= \text{dim}_{\bf C}H_{\bf C}+\frac{1}{2}\sum_{p,q}(p+q)^{2}(t^{p,q}_{H}-h^{p,q}_{H})w^{4}$$
et donc :
$$\emph{c}_{2}( {\xi}_{{\bf P}^2}(H_{\bf C},W_{.}^{\bullet },F^{\bullet },{\overline F}^{\bullet }))=\frac{1}{2}\sum_{p,q}(p+q)^{2}(h^{p,q}_{H}-t^{p,q}_{H})w^{4}.$$

\end{prop}

\begin{preuve}
On applique la proposition \ref{filtrscin} dans le cas où l'espace vectoriel trifiltré provient d'une structure de Hodge mixte avec $(V,F^{\bullet }_{0},F^{\bullet }_{1},F^{\bullet }_{2})=(H_{\bf C},W_{.}^{\bullet },F^{\bullet },{\overline F}^{\bullet })$. D'où  $  \text{dim}_{\bf C}Gr_{F_{2}^{\bullet }}^{q}Gr_{F_{1}^{\bullet }}^{p}V= \text{dim}_{\bf C}Gr_{{\overline F}_{}^{\bullet }}^{q}Gr_{{ F}_{}^{\bullet }}^{p}H_{\bf C}=t^{p,q}_{H}$ et $ \text{dim}_{\bf C}Gr_{F_{2}^{\bullet }}^{q}Gr_{F_{1}^{\bullet }}^{p}Gr_{F_{0}^{\bullet }}^{-p-q}V = \text{dim}_{\bf C}Gr_{{\overline F}^{\bullet }}^{q}Gr_{F_{}^{\bullet }}^{p}Gr^{W^{\bullet }_{.}}_{-p-q}H_{\bf C}$ ce qui donne la relation.
\end{preuve}
Ceci nous amène à poser pour définition du niveau de ${\bf R}$-scindement :   
\begin{definition}
Le niveau de ${\bf R}$-scindement d'une structure de Hodge mixte $(H_{\bf Q},W_{\bullet },F^{\bullet},{\overline F}^{\bullet})$ est le nombre entier : $$\alpha (H)=\emph{ch}_{2}({\xi}_{{\bf P}^2}(H_{\bf C},W_{.}^{\bullet },F^{\bullet },{\overline F}^{\bullet }))=\frac{1}{2}\sum_{p,q}(p+q)^{2}(t^{p,q}_{H}-h^{p,q}_{H}).$$
\end{definition}
Cette définition constitue une généralisation de la notion de structure de Hodge ${\bf R}$-scindée. En effet, si $(H_{\bf Q},W_{\bullet },F^{\bullet},{\overline F}^{\bullet})$ est une structure de Hodge ${\bf R}$-scindée alors $\alpha (H)=0$. Reciproquement supposons que $\alpha(H)=0$ pour une structure de Hodge mixte $H$. On verra par la suite que cela implique que pour tous $(p,q)$ $t_{H}^{p,q}=h^{p,q}_{H}$. On veut montrer que dans ce cas pour tous $(p,q)$ $I^{p,q}={\overline I}^{q,p}$. Supposons que ce ne soit pas le cas. Munissons ${\bf Z}^2$ de l'ordre lexicographique. Soit $(p_{0},q_{0})$ le plus grand des éléments $(p,q)$ tels que $I^{p,q}\neq I^{q,p}$ (rappelons que l'on a toujours égalité modulo $W_{p+q-2}$). Soit $I^{p_{0},q_{O}}_{0}$ le sous-espace vectoriel de $I^{p,q}$ de dimension maximale vérifiant l'égalité, donné par $I^{p_{0},q_{0}}=I^{p,q}\cap {\overline I}^{p,q}\subsetneq I^{p,q}$. On a alors $t^{p_{0},q_{0}}_{H}=dim_{{\bf C}}\,I^{p_{0},q_{O}}_{0} < dim_{\bf C}\,I^{p_{0},q_{O}}=h^{p_{0},q_{0}}_{H}$ ce qui contredit l'hypothèse. D'où :
\begin{center}
$ \alpha(H)=0$    $\Longleftrightarrow$  La structure de Hogde mixte $ H$ est  ${\bf R}$-scind\'ee.
\end{center}
\textit{Remarques :} $\bullet$ Notons $H'=(H_{\bf C},W_{\bullet },e^{-i.\delta}.F^{\bullet })$ la structure de Hodge mixte ${\bf R}$-scindée associée à la structure de Hodge mixte $H=(H_{\bf C},W_{\bullet },F^{\bullet })$. On vérifie que l'on a bien pour tous $(p,q)$ $h^{p,q}_{H'}=t^{p,q}_{H'}=dim_{\bf C}\,I^{p,q}$ et donc $\alpha(H_{\bf C},W_{\bullet },e^{-i.\delta}.F^{\bullet })=0$.

$\bullet$ Pour des raisons de dimensions, on a $\sum_{p,q}(h^{p,q}_{H}-t^{p,q}_{H})=0$. De plus $\text{c}_{1}({\xi}_{{\bf P}^2}(H_{\bf C},W_{.}^{\bullet },F^{\bullet },{\overline F}^{\bullet }))=\sum_{p,q}(p+q)(h^{p,q}_{H}-t^{p,q}_{H})=0$ comme anoncé pour un fibré de Rees associé à des filtrations opposées.

$\bullet $ Il para\^{\i}t possible que deux structures de Hodge mixtes $H$ et $H'$ aient les m\^emes nombres de Hodge, des niveaux de ${\bf R}$-scindement $\alpha (H)$ et $\alpha (H')$ égaux mais que leurs nombres $h^{p,q}_H$ et $t^{p,q}_H$ ne soient pas tous égaux.
 
\subsection{Comportement du niveau de ${\bf R}$-scindement par opérations sur les structures de Hodge mixtes}

Dans cette partie nous allons étudier les variations de $\alpha$ lors des différentes opérations qui peuvent \^etre faites sur des structures de Hodge mixtes héritées de la structure de catégorie abélienne de la catégorie des structures de Hogde mixtes. Morphismes entre structures de Hodge mixtes, quotients, sommes directes, produits tensoriels et extensions.

\begin{definition}
Soit $A$ une structure de Hodge mixte donn\'ee on note ${\mathcal E }_{A}$ le sous-ensemble de ${\bf Z} \times {\bf Z} $ formé par les couples $(p,q) \in {\bf Z}\times {\bf Z}$ tels que les nombres de Hodge $h^{p,q}_{A}$ sont non nuls si et ssi $(p,q) \in {\mathcal E_{A}}$ . Cet ensemble est appel\'e le type de $A$.
\end{definition}
\begin{definition}
La structure de Hodge $T \langle k \rangle$ est l'unique structure de Hodge de rang $1$, de type $(-k, -k)$ et de réseau entier $ (2 \pi i )^{k}{\bf Z}$.
\end{definition}

\textit{Remarque :} Soit $A$ une structure de Hodge mixte et $k \in {\bf Z}$, alors les types de $A$ et $A \otimes {T \langle k \rangle}$ sont liés par : ${\mathcal E }_{A \otimes {T \langle k \rangle}}= \{ (p-k,q-k)\in {\bf Z} \times {\bf Z} \vert (p,q) \in {\mathcal E }_{A} \} $.

Pour les filtrations on a les relations pour $p,q,n \in {\bf Z}^{3}$ :

\hspace{1cm}$\left\lbrace \begin{array}{l}
F^{p}(H \otimes {T \langle k \rangle})= F^{p+k}(H),\\   
{\overline F}^{q}(H \otimes {T \langle k \rangle})= {\overline F}^{q+k}(H),\\
W_{n}(H \otimes {T \langle k \rangle})= W_{n+2k}(H),    
\end{array}
\right.$ \hspace{0.5cm} et \hspace{0.5cm} $\left\lbrace \begin{array}{l}
t^{p,q}_{H \otimes {T \langle k \rangle}}=t^{p+k,q+k}_{H},\\
h^{p,q}_{H \otimes {T \langle k \rangle}}=h^{p+k,q+k}_{H }.\\   
\end{array}
\right.$

\begin{theorem}\label{elemalpha}
Pour $H$ et $H'$ deux structures de Hodge mixtes :

\item{(i)} pour tout $k \in {\bf Z}$, $\alpha(H \otimes T \langle k \rangle )=\alpha(H)$. 

\item{(ii)} $\alpha(H^{*})=\alpha(H)$ où $H^{*}=\emph{Hom}_{SHM}(H,T \langle 0 \rangle )$.

\item{(iii)} $\alpha ( H \oplus H')=\alpha(H) +\alpha (H')$.
 
\item{(iv)} $\alpha( H \otimes H')=\emph{dim}(H').\alpha(H)+\emph{dim}(H).\alpha(H')$. 
\end{theorem}

\begin{preuve}(i)\begin{eqnarray*}
\alpha(H \otimes T \langle k \rangle )&=&\frac{1}{2}\sum_{(p,q)\in {\cale}_{H \otimes T \langle k \rangle}}(p+q)^{2}(t^{p,q}_{H \otimes T \langle k \rangle}-h^{p,q}_{H \otimes T \langle k \rangle})\\
&=&\frac{1}{2}\sum_{(p+k,q+k)\in {\cale}_{H}}(p+q)^{2}(t^{p,q}_{H \otimes T \langle k \rangle}-h^{p,q}_{H \otimes T \langle k \rangle})\\
&=&\frac{1}{2}\sum_{(p,q)\in {\cale}_{H}}(p+q)^{2}(t^{p-k,q-k}_{H \otimes T \langle k \rangle}-h^{p-k,q-k}_{H \otimes T \langle k \rangle})\\
&=&\alpha(H),
\end{eqnarray*}
d'après la remarque précédente.\\
(ii) Il suffit d'écrire que ${\cale }_{H^*}=\{(-p,-q) \vert (p,q) \in {\cale}_{ H}$, $t^{p,q}_{H^*}=t^{-p,-q}_{H}$ et $h^{p,q}_{H^*}=h^{-p,-q}_{H}\}$, l'égalité en découle car tous les coefficients $p+q$ sont au carré.\\
(iii) ${\cale }_{H\oplus H'}={\cale}_{H} \cup {\cale}_{ H'}$. Les filtrations de Hodge et la filtration décroissante associée à la filtration par le poids de $H \oplus H'$ se déduisent des filtrations respectives de $H$ et $H'$ par somme directe. Les dimensions des quotients sont donc faciles à calculer, on en déduit : pour $p,q$ entiers, $h^{p,q}_{H \oplus H'}=h^{p,q}_{H}+t^{p,q}_{H'}$ et $t^{p,q}_{H \oplus H'}=t^{p,q}_{H}+t^{p,q}_{H'}$. En découle la formule voulue.\\
(iv) D'après le lemme \ref{dirsumxi} :\\
 ${\xi}_{{\bf P}^2}((H_{\bf C},W^{\bullet }_{.},F^{\bullet},{\overline F}^{\bullet}) \otimes (H_{\bf C} {}',W^{\bullet }_{.} {}',F^{\bullet} {}',{\overline F}^{\bullet} {} ')) \cong {\xi}_{{\bf P}^2}((H_{\bf C} {},W^{\bullet }_{.} {},F^{\bullet} {},{\overline F}^{\bullet} {} )) \otimes {\xi}_{{\bf P}^2}((H_{\bf C} {}',W^{\bullet }_{.} {}',F^{\bullet} {}',{\overline F}^{\bullet} {} ')) $.\\
Donc :

$\text{ch}{\xi}_{{\bf P}^2}((H_{\bf C},W^{\bullet }_{.},F^{\bullet},{\overline F}^{\bullet}) \otimes (H_{\bf C} {}',W^{\bullet }_{.} {}',F^{\bullet} {}',{\overline F}^{\bullet} {} '))=$

\hspace{4cm} $\text{ch}({\xi}_{{\bf P}^2}((H_{\bf C} {},W^{\bullet }_{.} {},F^{\bullet} {},{\overline F}^{\bullet} {} ))).\text{ch}({\xi}_{{\bf P}^2}((H_{\bf C} {}',W^{\bullet }_{.} {}',F^{\bullet} {}',{\overline F}^{\bullet} {} '))) $.\\
Or les trois filtrations constituant une structure de Hodge étant mixte étant opposées, on a : 
$\text{ch}_{1}({\xi}_{{\bf P}^2}((H_{\bf C},W^{\bullet }_{.},F^{\bullet},{\overline F}^{\bullet}))=0$ et de m\^eme $\text{ch}_{1}({\xi}_{{\bf P}^2}((H_{\bf C}{} ',W^{\bullet }_{.}{} ',F^{\bullet} {}',{\overline F}^{\bullet}{} '))=0$. Ainsi :

$\text{c}_{2}({\xi}_{{\bf P}^2}((H_{\bf C},W^{\bullet }_{.},F^{\bullet},{\overline F}^{\bullet}) \otimes (H_{\bf C} {}',W^{\bullet }_{.} {}',F^{\bullet} {}',{\overline F}^{\bullet} {} ')))=$

\hspace{1cm} $\text{dim}_{\bf C}(H_{\bf C}{}').\text{c}_{2}({\xi}_{{\bf P}^2}((H_{\bf C} {},W^{\bullet }_{.} {},F^{\bullet} {},{\overline F}^{\bullet} {} )))+\text{dim}_{\bf C}(H_{\bf C}).\text{c}_{2}({\xi}_{{\bf P}^2}((H_{\bf C} {}',W^{\bullet }_{.} {}',F^{\bullet} {}',{\overline F}^{\bullet} {} '))) $.\\
Qui est l'égalité cherchée.

\end{preuve}

Soient $A=(A_{\bf Z},W_{\bullet }^{A},F^{\bullet}_{A},{\overline F}^{\bullet}_{A})$ et $B=(B_{\bf Z},W_{\bullet }^{B},F^{\bullet}_{B},{\overline F}^{\bullet}_{B})$ deux structures de Hodge mixtes. Nous pouvons alors définir  suivant \cite{mor} ou \cite{car} le groupe d'extension de $B$ par $A$ noté $\text{Ext}(B,A)$.
C'est à dire l'ensemble des classes d'équivalence à congruence prés des suites exactes de structures de Hodge mixtes :
$$
\xymatrix{
 0 \ar[r] &  A        \ar[r]^{i}  & H         \ar[r]^{\pi}     &   B    \ar[r]   & 0  
}
$$

Pour une structure de Hodge mixte $H$ remplissant une telle suite exacte, on écrira : $H \in  \text{Ext(B,A)}$. On a bien s\^ur :
 ${\mathcal E }_{H}={\mathcal E }_{A} \cup {\mathcal E }_{B}$.
Nous noterons pour tous $p,q \in {\bf Z}$ $t^{p,q}_{H}$ (resp. $t^{p,q}_{A}$, $t^{p,q}_{B}$) les entiers suivants $\text{dim}_{\bf C}Gr_{{\overline F}_{H}^{\bullet }}^{q}Gr_{{ F}_{H}^{\bullet }}^{p}H_{\bf C}$ (resp. $\text{dim}_{\bf C}Gr_{{\overline F}_{A}^{\bullet }}^{q}Gr_{F_{A}^{\bullet }}^{p}A_{\bf C}$, $\text{dim}_{\bf C}Gr_{{\overline F}_{B}^{\bullet }}^{q}Gr_{F_{B}^{\bullet }}^{p}B_{\bf C}$). On utilisera le m\^eme type de  notation pour les nombres de Hodge $h^{p,q}_{H}$,$h^{p,q}_{A}$ et $h^{p,q}_{B}$.\\
Pour $H \in \text{Ext}(B,A)$ nous avons les égalités suivantes pour les nombres de Hodge :
\begin{center}
Pour tous $p,q \in {\bf Z}$, $h^{p,q}_{H}=h^{p,q}_{A}+h^{p,q}_{B}$
\end{center}    
Une telle égalité pour les nombres $t^{p,q}$ n'est pas vraie en général (cf exemples plus bas). Ecrivons le diagramme :



\scalebox{0.70}[1]{
$$
\xymatrix{
{}  &    0   \ar[d] &   0   \ar[d] & 0  \ar[d] &  {}\\
 0  \ar[r] &  {\xi}_{\widetilde{{\bf P}^2}}(A_{\bf C},{W^{A}_{. }}^{\bullet },F_{A}^{\bullet },{\overline F}_{A}^{\bullet },Triv^{\bullet })       \ar[r]  \ar[d] & {\xi}_{\widetilde{{\bf P}^2}}(H_{\bf C},{W^{H}_{.}}^{\bullet },F_{H}^{\bullet },{\overline F}_{H}^{\bullet },Triv^{\bullet })         \ar[r] \ar[d]     &   {\xi}_{\widetilde{{\bf P}^2}}(B_{\bf C},{W^{B}_{.}}^{\bullet },F_{B}^{\bullet },{\overline F}_{B}^{\bullet },Triv^{\bullet })    \ar[r] \ar[d]  & 0  \\
          {}  &  e^{*}{\xi}_{{{\bf P}^2}}(A_{\bf C},{W^{A}_{.}}^{\bullet },F_{A}^{\bullet },{\overline F}_{A}^{\bullet })  \ar[d]      &   e^{*}{\xi}_{{{\bf P}^2}}(H_{\bf C},{W^{H}_{.}}^{\bullet },F_{H}^{\bullet },{\overline F}_{H}^{\bullet })    \ar[d]      &  e^{*}{\xi}_{{{\bf P}^2}}(B_{\bf C},{W^{B}_{.}}^{\bullet },F_{B}^{\bullet },{\overline F}_{B}^{\bullet })      \ar[d]   &   {}\\
{}  &  {\calf}_{A}   \ar[d]      &   {\calf}_{H}    \ar[d]      &       {\calf}_{B}  \ar[d]   &   {}\\
{}  & 0         &  0          &      0   &   {}
}
$$
}

   







où, d'après la section précédente, les colonnes sont des suites exactes associées à chacunes des structures de Hodge mixtes $A$, $H$ et $B$. La première ligne est une suite exacte car le foncteur ${\xi}_{\widetilde{{\bf P}^2}}(.,.,.,.,Triv^{\bullet })$ est un foncteur exact de la catégorie des structures de Hodge mixtes vers la catégorie des fibrés vectoriels sur ${\widetilde{{\bf P}^2}}$. On obtient donc la relation entre les niveaux de ${\bf R}$-scindement des structures de Hodge mixtes :
$$  \alpha (H)= \alpha (A) + \alpha (B)+ \frac{1}{2}\sum_{p,q}\,(p+q)^{2}(t^{p,q}_{A}+t^{p,q}_{B}-t^{p,q}_{H}).$$
Cette relation est importante dans la démonstration du fait que l'invariant $\alpha$ est sur-additif par extensions :

\begin{theorem}\label{inegext}
Soient $A$ et $B$ deux structures de Hodge mixtes et $H \in \emph{Ext}(B,A)$, alors :
$$ \alpha(H) \geq \alpha(A)+\alpha(B)$$

\end{theorem} 

 Notation : pour une structure de Hodge mixte $\bullet $, on défini $f^{p,q}_{\bullet}=\text{dim}_{\bf C}F^{p}_{\bullet}\cap{\overline F}^{q}_{\bullet}$. La démonstration du théorème utilise le lemme suivant ainsi que son corollaire :
  
\begin{lemma}
Soient $V_1$ et $V_2$ deux espaces vectoriels, $V=V_{1}\oplus V_{2}$ et $\pi$ la projection de $V$ sur $V_2$, $i$ l'injection de $V_1$ dans $V$. Soient $W_1$, $W_{1}'$ (resp. $W_2$, $W_{2}'$) des sous-espaces vectoriels de $V_1$ (resp. $ V_{2}$). Si $W$ et $W'$ sont des sous-espaces vectoriels tels que $i(W_{1})=W \cap V_{1},i(W_{1}')=W' \cap V_{1}$ et $\pi (W)=W_{2},\pi (W')=W_{2}'$, alors :

$  \emph{dim}_{\bf C}(W_{1} \cap W_{1}')+ \emph{dim}_{\bf C}(W_{2} \cap W_{2}')-\min(\emph{dim}_{\bf C}(W_{2} \cap W_{2}'),\emph{dim}_{\bf C}(V_{1}))$
\begin{eqnarray*}
 \hspace{4cm}  &   \leq   &   \emph{dim}_{\bf C}( W\cap W')\\
\hspace{4cm} & \leq & \emph{dim}_{\bf C}(W_{1} \cap W_{1}')+\emph{dim}_{\bf C}(W_{2} \cap W_{2}'). 
\end{eqnarray*}

\end{lemma}

\begin{corollaire}
Soient $A$ et $B$ deux structures de Hodge mixtes et $H \in \emph{Ext}(B,A)$, alors :
$$ \forall (p,q) \in {\bf Z} \times {\bf Z} \,\,\,\, f^{p,q}_{H}-f^{p,q}_{A}-f^{p,q}_{B} \leq 0$$
 
\end{corollaire}

\begin{preuve}(du corollaire)
D'après la construction de $H \in \text{Ext}(B,A)$ (c'est à dire une classe de suite exacte  $0 \longrightarrow A \stackrel{i}{\longrightarrow} H \stackrel{ \pi}{\longrightarrow} B \longrightarrow 0 $) dans \cite{pen} suivant \cite{mor},\cite{car}, pour tout $(p,q) \in {\bf Z} \times {\bf Z}$, on est exactement dans le cadre du lemme précédent avec : $H_{\bf C}=A_{\bf C}\oplus B_{\bf C}$, $W_1=F^{p}_{A}$, $W_{1}'={\overline F}^{q}_{A}$, $W_1=F^{p}_{B}$, $W_{1}'={\overline F}^{q}_{B}$, $W=F^{p}_{H}$, $W'={\overline F}^{q}_{H}$. Par construction de la filtration de Hodge et sa filtration opposée sur l'extension, on a
$i(F^{p}_{A_{\bf C}})=A_{\bf C} \cap F^{p}_{H_{\bf C}}$, $i({\overline F}^{q}_{A_{\bf C}})=A_{\bf C} \cap {\overline F}^{q}_{H_{\bf C}}$ et $\pi(F^{p}_{H_{\bf C}})=F^{p}_{B_{\bf C}}$, $\pi({\overline F}^{q}_{H_{\bf C}})= {\overline F}^{q}_{B_{\bf C}}$ par stricte compatibilité des morphismes de structures de Hodge.
\end{preuve}

\begin{preuve}(du lemme)
Ecrivons sous forme matricielle les coordonnées des sous espaces vectoriels $W$ et $W'$ dans $V=V_{1}\oplus V_{2}$ i.e les représentations matricielles des points $W$ et $W'$ dans les grassmanniennes $G(V,\text{dim}_{\bf C}W)$ et $G(V,\text{dim}_{\bf C}W')$. On note $M_{i}Sev$ la matrice représentant $Sev \subset V_{i}$ dans $G(V_{i},\text{dim}_{\bf C}Sev)$ pour $i \in \{\Box ,1,2\}$. La base de $V$ prise pour la représentation matricielle est la réunion des bases de $V_{1}$ et $V_{2}$. Alors :\\
$MW=\left(
    \begin{array}{ccc}
    M_{1}W_{1}   &   \vert    &      0   \\
    A_{1}     &    \vert       &     M_{2}W_{2}

\end{array}
\right)$    et $MW'=\left(
    \begin{array}{ccc}
    M_{1}W_{1}'   &   \vert    &      0   \\
    A_{1}'     &    \vert       &     M_{2}W_{2}'

\end{array}
\right)$.\\
où $A_{1}$ et $A_{2}$ sont des matrices quelconques de dimension $\text{dim}_{\bf C}W_{2} \times \text{dim}_{\bf C}V_{1} $ et $\text{dim}_{\bf C}W_{2}' \times \text{dim}_{\bf C}V_{1}$ respectivement (on vérifie que ce sont bien des sous-espaces vectoriels tels que $i(Sev_{1})=Sev \cap V_{1}$ et $\pi(Sev)=Sev_{2}$). Comme pour deux sous-espaces vectoriels on a l'égalité sur les dimensions : $\text{dim}(Sev_{1}\cap Sev_{2})+\text{dim}(Vect(Sev_{1},Sev_{2}))=\text{dim}(Sev_{1})+dim(Sev_{2})$, conna\^{\i}tre la dimension de $Vect(W,W')$ nous donne la dimension cherchée. Il s'agit donc de trouver le maximum des dimensions des matrices extraites de déterminant non nul de la matrice : 

\hspace{2cm}$MW=\left(
    \begin{array}{ccc}
    M_{1}W_{1}   &   \vert    &      0   \\
    A_{1}     &    \vert       &     M_{2}W_{2}\\
    --------&   \vert  & --------\\
    M_{1}W_{1}'   &   \vert    &      0   \\
    A_{1}'     &    \vert       &     M_{2}W_{2}'

\end{array}
\right)$\\
Une telle matrice contient des lignes et colonnes obtenues à partir des matrices extraites de déterminant non nul de dimension maximale pour les sous-espaces engendrés $Vect( W_{1},  W_{1}')$ et $Vect( W_{2},  W_{2}')$ i.e les lignes et colonnes complétées à partir de telle matrices extraites dans :

\hspace{2cm} $MW'=\left(
    \begin{array}{c}
    M_{1}W_{1} \\
    M_{1}W_{1}' 
\end{array}
\right)$ et $MW'=\left(
    \begin{array}{c}
    M_{2}W_{2} \\
    M_{2}W_{2}' 
\end{array}
\right)$.\\
On peut donc trouver une matrice extraite de déterminant non nul de dimension au moins égale à $ \text{dim}\,Vect( W_{1},  W_{1}')+\text{dim}\,Vect( W_{2},  W_{2}')$. On peut ajouter à cela au plus $\text{min}(\text{dim}\,( W_{2} \cap  W_{2}'), \text{dim} \,V_{1})$ lignes et colonnes pour obtenir une matrice de déterminant non nul. D'où l'inégalité :
\begin{center}
$ \text{dim}(Vect(W_{1},  W_{1}'))+\text{dim}(Vect(W_{2},  W_{2}'))   \leq  \text{dim}(Vect(W_{},  W_{}')) $\\
$  \leq  \text{dim}(Vect(W_{1},  W_{1}'))+\text{dim}(Vect(W_{2},  W_{2}'))+ \text{min}(\text{dim}\,( W_{2} \cap  W_{2}'), \text{dim} \,V_{1})$
\end{center}
En découle l'inégalité voulue en passant aux codimensions. 


\end{preuve}

\begin{preuve}(du théorème)
 Comme ${\mathcal E }_{H}={\mathcal E }_{A} \cup {\mathcal E }_{B}$ est un sous-ensemble fini de ${\bf Z} \times {\bf Z}$, on peut trouver $k \in {\bf Z}$ et $N \in {\bf Z}$ tels que ${\mathcal E }_{H \otimes {T \langle k \rangle }}={\mathcal E }_{A \otimes {T \langle k \rangle }} \cup {\mathcal E }_{B\otimes {T \langle k \rangle }} \subset [0,N] \times [0,N]$. D'après le théorème \ref{elemalpha} $(i)$, l'égalité n'est pas modifiée par tensorisation des structures de Hodge $A$,$B$ et $H$ par une structure de Hodge de Tate. Tensoriser les structures de Hodge dans une extension par une structure de Hodge de Tate ${T \langle k \rangle }$ induit un isomorphisme entre $\text{Ext}(B,A)$ et $\text{Ext}(B\otimes {T \langle k \rangle },A\otimes {T \langle k \rangle })$. Quitte à remplacer les structures de Hodge $A,B$ et $H$ par leurs tensorisées par ${T \langle k \rangle }$, on peut ainsi supposer que  ${\mathcal E }_{H } \subset [0,N] \times [0,N]$. La proposition se ramène donc à montrer que : 
$$ \alpha(H) - \alpha(A)-\alpha(B) =\frac{1}{2} \sum_{(p,q) \in {[0,N] \times [0,N]}} (p+q)^{2}(t^{p,q}_{A}+t^{p,q}_{B}-t^{p,q}_{H}) \leq 0$$
Pour une structure de Hodge mixte $\bullet$, on a les suites exactes suivantes :

$$
\xymatrix{
{}   &   0  \ar[d]    &    0   \ar[d]   &   {}    &   {}\\
{}   &   F^{p+1}_{\bullet}\cap{\overline F}^{q+1}_{\bullet}  \ar[d]    &    F^{p}_{\bullet}\cap{\overline F}^{q+1}_{\bullet}   \ar[d]   &   {}    &   {}\\
  {}   &   F^{p+1}_{\bullet}\cap{\overline F}^{q}_{\bullet}  \ar[d]    &    F^{p}_{\bullet}\cap{\overline F}^{q}_{\bullet}   \ar[d]   &   {}    &   {}\\
 0  \ar[r]  &   F^{p+1}_{\bullet}Gr_{{\overline {F_{\bullet}}}}^{q}   \ar[r] \ar[d] &  F^{p}_{\bullet}Gr_{{\overline {F_{\bullet}}}}^{q}   \ar[r] \ar[d]   &  Gr_{{ {F_{\bullet}}}}^{p}Gr_{{\overline {F_{\bullet}}}}^{q}  \ar[r]  &   0\\
{}    &   0    &    0    &  {}   &   {}
}
$$
Et donc la relation sur les dimensions :$$ t^{p,q}_{\bullet}=f^{p,q}_{\bullet}-f^{p+1,q}_{\bullet}-f^{p,q+1}_{\bullet}+f^{p+1,q+1}_{\bullet}.$$
Ainsi :
$$\alpha^{-}( {\bullet})= \frac{1}{2} \sum_{(p,q) \in {[0,N] \times [0,N]}} (p+q)^{2}\,t^{p,q}_{\bullet }=\frac{1}{2}\sum_{(p,q) \in {[0,N] \times [0,N]}} (p+q)^{2}(f^{p,q}_{\bullet}-f^{p+1,q}_{\bullet}-f^{p,q+1}_{\bullet}+f^{p+1,q+1}_{\bullet}).$$ 
On effectue un changement d'indice sur les sommes, comme pour tout $p$, $f^{p,N+1}_{\bullet}=0$, pour tout $q$, $f^{N+1,q}_{\bullet}=0$ :
\begin{center}
$\alpha^{-}( {\bullet})= \frac{1}{2} \sum_{(p,q) \in {[1,N] \times [1,N]}} \Bigl[ (p+q)^{2}-(p+q-1)^{2}-(p+q-1)^{2}+(p+q-2)^{2}              \Bigr] f^{p,q}_{\bullet}+\frac{1}{2} \sum_{q \in [0,N]}\Bigl[ (0+q)^{2}-(0+q-1)^{2}            \Bigr] f^{0,q}_{\bullet}$\\

$+\frac{1}{2} \sum_{p \in [1,N]}\Bigl[ (p+0)^{2}-(p-1+0)^{2}            \Bigr] f^{p,0}_{\bullet} $.
\end{center}
$\alpha^{-}( {\bullet})= \frac{1}{2} \sum_{(p,q) \in {[1,N] \times [1,N]}} ( 2pq+4) f^{p,q}_{\bullet}+\frac{1}{2} \sum_{q \in [0,N]} ( 2q-1) f^{0,q}_{\bullet}+\frac{1}{2} \sum_{p \in [1,N]}(2p-1) f^{p,0}_{\bullet} $,\\
que nous écrirons, afin que les coefficients des $f^{p,q}_{\bullet}$ sous les $\sum$ soient tous positifs, sous la forme :\\
$\alpha^{-}( {\bullet})= \frac{1}{2} \sum_{(p,q) \in {[1,N] \times [1,N]}} ( 2pq+4) f^{p,q}_{\bullet}+\frac{1}{2} \sum_{q \in [1,N]} ( 2q-1) f^{0,q}_{\bullet}+\frac{1}{2} \sum_{p \in [1,N]}(2p-1) f^{p,0}_{\bullet}- \frac{1}{2} f^{0,0}_{\bullet}$\\
Ainsi :\\
$\alpha( H)-\alpha(A)-\alpha(B)=\alpha^{-}( H)-\alpha^{-}(A)-\alpha^{-}(B)= \frac{1}{2} \sum_{(p,q) \in {[1,N] \times [1,N]}} ( 2pq+4) \Bigl[ f^{p,q}_{H}-f^{p,q}_{A}-f^{p,q}_{B} \Big]+\frac{1}{2}\sum_{q \in [1,N]} ( 2q-1) \Bigl[ f^{0,q}_{H}-f^{0,q}_{A}-f^{0,q}_{B} \Big]+\frac{1}{2} \sum_{p \in [1,N]}(2p-1) \Bigl[ f^{p,0}_{H}-f^{p,0}_{A}-f^{p,0}_{B} \Big]- \frac{1}{2} \Bigl[ f^{0,0}_{H}-f^{0,0}_{A}-f^{0,0}_{B} \Big]$.\\
$f^{0,0}_{H}-f^{0,0}_{A}-f^{0,0}_{B}=\text{dim}_{\bf C}H-\text{dim}_{\bf C}A-\text{dim}_{\bf C}B=0$, donc tous les coefficients des $f^{p,q}_{\bullet}$ dans les $\sum$ étant positifs, d'après le lemme précédent, pour tous $(p,q) \in [0,N] \times [0,N]$, $f^{p,q}_{H}-f^{p,q}_{A}-f^{p,q}_{B} \leq 0$. D'où le résultat.


\end{preuve}

\begin{corollaire}
Soit $H$ une structure de Hodge mixte, alors :    \,\,\,\,       $\alpha(H) \geq 0$.

\end{corollaire}

\begin{preuve}
Toute structure de Hodge mixte peut \^etre décrite comme extension successive de structures de Hodge pures. En effet, si $H$ est une structure de Hodge dont les poids varient entre $0$ et $2n$ par exemple (on peut s'y ramener en tensorisant par une structure de Hodge de Tate, ce qui ne change pas la valeur de $\alpha$ d'après la proposition \ref{elemalpha}), on a $H \in \text{Ext}(Gr^{W}_{2n}H,W_{2n-1}H)$. Ainsi par la proposition précédente $\alpha(H) \geq \alpha(W_{2n-1}H) + \alpha(Gr^{W}_{2n}H)=\alpha(W_{2n-1}H)$ car $\alpha(Gr^{W}_{2n}H)=0$, $Gr^{W}_{2n}H$ étant une structure de Hodge pure de poids $2n$. De m\^eme $W_{2n-1}H \in \text{Ext}(Gr^{W}_{2n-1}H,W_{2n-2}H)$ et donc $\alpha(H) \geq \alpha(W_{2n-1}H) \geq \alpha(W_{2n-2}H)$. En poursuivant ainsi :
$$  \alpha(H) = \alpha(W_{2n}H) \geq \alpha(W_{2n-1}H) \geq \alpha(W_{2n-2}H) \geq ...\geq \alpha(W_{0}H)=0,$$
car $W_{0}H$ est une structure de Hoge pure.
\end{preuve}

\textit{Remarques :} $\bullet$ En appliquant le critère de Drézet et de Le Potier, on voit que les fibrés vectoriels $\xi_{{\bf P}^2}(H^{k}(X,{\bf C}),W^{\bullet }_{.},F^{\bullet },{\overline F}^{\bullet })$ sont semistables. En effet, on est dans le cas où $r>1$,$c_1=0$ et $c_2 \leq 0$.\\ 
$\bullet$ Soient $A$ et $B$ deux structures de Hodge mixtes, il existe  $m(A,B) \in {\bf Z}$ ($m=m(h_{A}^{.,.},h_{B}^{.,.},t_{A}^{.,.},t_{B}^{.,.}$)) tel que pour tout $H \in \text{Ext}(B,A)$ : $\alpha(H) \in [\alpha(A)+\alpha(B),m(A,B)]$. On a m\^eme plus : pour tout $p \in [\alpha(A)+\alpha(B),m(A,B)]$ il existe $H \in \text{Ext}(B,A)$ tel que $\alpha(H)=p$. On en déduit qu'étant donnés des nombres de Hodge $h^{p,q}$, il existe $m(h^{.,.})$ tel que si $H$ est une structure de Hodge mixte qui a ces nombres de Hodge alors $\alpha(H) \in [0,m(h^{.,.})]$.


\subsection{Exemples géométriques}

D'après Deligne \cite{del2}, \cite{del3}, les groupes de cohomologie des variétés algébriques (schémas séparés de type fini sur ${\bf C}$) sont munis d'une struture de Hodge mixte. Soit $X$ une telle variété, on peut lui associer pour tout $k \in {\bf Z}$ des entiers :
$$\alpha_{k}(X)=\alpha(H^{k}=(H^{k}(X,{\bf C}),W^{\bullet }_{.},F^{\bullet },{\overline F}^{\bullet }))$$ 
Le théorème précédent permet de préciser le comportement de ces entiers vis-à-vis de certaines morphismes applications entre variétés algébriques :

\begin{corollaire}
\begin{itemize}\parindent=2cm
\item {(i)} Soit $X$ une variété algébrique, $X'$ sa normalisée et $\overline X$ sa complétée, alors pour tout $k \in {\bf Z}$ :
\begin{center}
$\alpha_{k}(X) \geq \alpha_{k}(X')$ et $\alpha_{k}(X) \geq \alpha_{k}({\overline X})$.
\end{center} 
\item {(ii)} Plus généralement, si $f:X \rightarrow Y$ est une application qui induit un morphisme injectif (resp. surjectif) sur la cohomologie $f^{*}:H^{k}(Y,{\bf C}) \rightarrow H^{k}(X,{\bf C})$ alors, pour tout $k \in {\bf Z}$ :
\begin{center}
$\alpha_{k}(X) \geq \alpha_{k}(Y)$ (resp. $\alpha_{k}(Y) \geq \alpha_{k}({ X})$).
\end{center} 
\par\end{itemize}
\end{corollaire}

\begin{preuve}
(i) se déduit de (ii) à l'aide de \cite{del3} prop (8.2.6) : l'application de la normalisée $X'$ vers $X$ induit un morphisme surjectif en cohomologie, l'application de $X$ vers sa complétée induit un morphisme injectif en cohomologie. La catégorie des structures de Hodge mixtes étant abélienne, on peut écrire des suites exactes à partir des injections et surjections de structures de Hodge mixtes et donc pour prouver (ii), on écrit des extensions et on utilise le théorème \ref{inegext}.
\end{preuve}

\begin{corollaire}
Soient $X$ et $Y$ deux variétés algébriques, alors :
\begin{center}
$\alpha_{k}(X \times Y)= \sum_{i=0}^{i=k}\Bigl[ \ \emph{dim}_{\bf C}(H^{k-i}(Y)).\alpha_{i}(X)+\emph{dim}_{\bf C}(H^{i}(X)).\alpha_{k-i}(Y) \Bigr]$.
\end{center}
\end{corollaire}

\begin{preuve}
D'après \cite{del3} proposition (8.2.10), les isomorphisme de K\"unneth $H^{.}(X\times Y,{\bf C}) \cong  H^{.}(X,{\bf C}) \otimes H^{.}(Y,{\bf C})$ sont des isomorphismes de structures de Hodge. On applique le $(iii)$ du théorème \ref{elemalpha} à la somme donnée par la formule de K\"unneth, puis le $(iv)$ à chacun des facteurs.
\end{preuve}

\subsubsection{Calculs de $\alpha$ pour les courbes}

On cherche ici à calculer la matrice des périodes de courbes algébriques pour en déduire la valeur des entiers $t^{p,q}$ associés au premier groupe de cohomologie puis calculer la valeur de $\alpha$ de la structure de Hodge mixte associée.

Soit $X$ une courbe alg\'ebrique sur ${\bf C}$, $X'$ la normalis\'ee de $X$ et $r: X' \rightarrow X$ le morphisme qui s'en d\'eduit.
Soit ${\overline X}'$ la courbe projective non singuli\`ere dont $X'$ est un ouvert dense et $\overline X$ la courbe d\'eduite de
${\overline X}'$ en contractant chacun des $r^{-1}(s)$ pour $s \in X$.
Notons $\overline r$ le morphisme ${\overline X}' \rightarrow {\overline X}$, et $j$,$j'$ les inclusions respectives de $X$ dans
$\overline X$ et  $X'$ dans ${\overline X}'$.
Soit $S$ l'ensemble fini ${\overline X}-X$. D'où le carré cartésien suivant :
$$
\xymatrix{
     X' \ar@{^{(}->}[r]^{j'} \ar[d]_r  & {\overline X}' \ar[d]^{\overline r} \\
     X \ar@{^{(}->}[r]_{j} & {\overline X}
   }
$$

\begin{prop}\cite{del2}  La cohomologie de $X$ est donn\'ee par:

$H^{1}(X,{\bf C})={\bf H}^{1}({\overline X},[ {\mathcal O}_{\overline X} \stackrel{d}{ \rightarrow } r_{*} {\Omega }_{{\overline
X}'}^{1}(\log S)])$\\
De plus, la filtration par le poids est donnée par :

$W^{1}(H^{1}(X,{\bf C}))=\text{Im}(H^{1}({\overline X},{\bf C}) \rightarrow H^{1}(X,{\bf C})$

$W^{0}(H^{1}(X,{\bf C}))=\text{Ker}(H^{1}({\overline X},{\bf C}) \rightarrow H^{1}({\overline X}',{\bf C})$\\
La suite spectrale d\'efinie par la filtration b\^ete de $[ {\mathcal O}_{\overline X} \stackrel{d}{\rightarrow} r_{*}
{\Omega}_{{\overline X}'}^{1}(\log S)]$ d\'eg\'en\`ere  en $E_{2}^{.,.}$ et aboutit \`a la filtration de Hodge $F^{\bullet }$ de $H^{\bullet }(X,{\bf C})$.
\end{prop}
Calculons l'hypercohomologie du complexe $[ {\mathcal O}_{\overline X} \stackrel{d}{\rightarrow} r_{*}{\Omega}_{{\overline X}'}^{1}(\log S)]$. Pour un complexe de faisceaux $({\mathcal K^{\bullet }},d)$ sur une vari\'et\'e $X$ muni d'un recouvrement ${\bf \calu}=(\calu_{i})_{i \in I}$, l'hypercohomologie ${\bf H}^{\bullet }(X,{\mathcal K}^{\bullet })$ du complexe est d\'efinie comme \'etant la cohomologie du complexe total $ C^{(p,q)}=C^{p}({\bf \calu},\delta ,d )$ o\`u $d $ est l'op\'erateur du complexe et $\delta $ celui de $\check{C}$ech associ\'e au recouvrement. ${\bf H}^{\bullet }(X,{\mathcal K}^{\bullet })$ est l'aboutissement des suites spectrales d\'efinies par les filtrations \'evidentes de $C^{(p,q)}$ (cf \cite{grihar}). Le complexe double associ\'e \`a $[ {\mathcal O}_{\overline X} \stackrel{d}{\rightarrow} r_{*}{\Omega}_{{\overline
X}'}^{1}(\log S)]$ est donc ici, relativement à ${\bf \calu}$ :


$$
\xymatrix{
{C^{0}({\coprod_{i} \calu_{i}} ,{\mathcal O}_{\overline X})} \ar[r]^{\delta} \ar[d]_{d}     &    {C^{1}({\coprod_{i,j} \calu_{i} \cap \calu_{j}},{\mathcal O}_{\overline X})}   \ar[d]_{d}\\
C^{0}({\coprod_{i} \calu_{i}} ,r_{*}{\Omega}_{{\overline X}'}^{1}(\log S)) \ar[r]_-{\delta}  & C^{1}({\coprod_{i,j} \calu_{i} \cap \calu_{j}},r_{*}{\Omega}_{{\overline X}'}^{1}(\log S))
}
$$

Pour voir dans $H^{1}(X,{\bf C})$ les différents sous-espaces associés à la filtration par le poids $W_{0} ,W_{1} ,W_{2} $, on utilise la proposition pr\'ec\'edente qui se traduit par : $W_{1}H^{1}({X},{\bf C})={\bf H}^{1}({\overline X},[ {\mathcal O}_{\overline X} \stackrel{d}{\rightarrow} r_{*}
{\Omega}_{{\overline X} '}^{1}])$ et comme annonc\'e : $W_{2}H^{1}(X,{\bf C})={\bf H}^{1}({\overline X},[{\mathcal O}_{\overline X} \stackrel{d}{\rightarrow} r_{*} {\Omega}_{{\overline
X}'}^{1}(\log S)])$.

\subsubsection{Exemple de ${\bf P}^1$ avec quatres points marqués}

Appliquons la décomposition précédente au calcul de la cohomologie ${ H}^1(X,{\bf C})$ de la courbe $X$ obtenue à partir de ${\bf P}^1$ avec quatres points distincts distingués $\{p_{1},p_{2},P_{1},Q_{1} \}$(dans l'ordre) dont on a enlevé les deux premiers $p_{1}$ et $p_{2}$ et recollé les deux autres $P_{1}$,et $Q_{1}$. () justification du recollement.\\  
Pour effectuer le calcul, utilisons le recouvrement standard de  ${\bf P}^1$ par les cartes affines de coordonnées $\calu ={\bf P}^1-\{\infty \}$ munie de la coordonnée $u$ et $\calv ={\bf P}^1-\{0\}$ munie de la coordonnée $v$. On notera $\calu \calv$ l'ouvert intersection des deux cartes. Nous commencerons par calculer la cohomologie de courbe complétée $\overline X$.

Soit  $f \in C^{1}({\mathcal {UV}},{\mathcal O}_{\overline X})$ (on choisira $u=v^{-1}$ comme coordonn\'ee d'écriture sur ${\mathcal {UV}}$). Alors $f$ et $df$ peuvent s'\'ecrirent : $f(u)= \sum_{- \infty}^{\infty} a_{n}u^{n}$ et $df(u)= \sum_{- \infty}^{\infty} na_{n}u^{n-1}$.

Prenons pour $g \in C^{0}({\mathcal U} \coprod {\mathcal V},r_{*}{\Omega}_{{\overline X}'}^{1})$, $g_{0}(u)$ \'egal \`a la partie de
$dg(u)$ \`a exposants positifs en $u$ et $g_{1}(v)$ la partie de $dg(u)$ correspondant aux exposants n\'egatifs en $u$. Alors, on a, pour tout $u \in \calu \calv$ : $\delta (g_{0},g_{1})(u)=dg(u)$. Donc tout \'el\'ement de $C^{1}({\mathcal {UV}},{\mathcal O}_{\overline X})$ voit son image par $d$ dans $C^{1}({\mathcal {UV}},r_{*}{\Omega}_{{\overline X}'}^{1})$ annul\'ee par l'image par $\delta$ d'un \'el\'ement de $C^{0}({\mathcal U} \coprod {\mathcal V},r_{*}{\Omega}_{{\overline X}'}^{1})$, donc l'hypercohomologie ${\bf H}^{1}({\overline X},[ {\mathcal O}_{\overline X} \stackrel{d}{\rightarrow} r_{*}{\Omega}_{{\overline X} '}^{1}])$ c'est à dire la cohomologie du complexe :

\hspace{1cm}$
\xymatrix{
C^{0}({\mathcal U} \coprod {\mathcal V},{\mathcal O}_{\overline X})  \ar[r]^-{(\delta,d)}& {C^{1}({\mathcal {UV}},{\mathcal O}_{\overline X})} \oplus C^{0}({\mathcal U} \coprod {\mathcal V},r_{*}{\Omega}_{{\overline X}'}^{1})  \ar[r]^-{d+\delta} & C^{1}({\mathcal {UV}},r_{*}{\Omega}_{{\overline X}'}^{1})
}
$\\
est d\'etermin\'e par l'image par $\delta$ de $C^{0}({\mathcal U} \coprod {\mathcal V},{\mathcal O}_{\overline X})$ dans $C^{1}({\mathcal {UV}},{\mathcal
O}_{\overline X})$. Cherchons donc quelle sont les images des  \'el\'ements $(h_{0},h_{1}) \in C^{0}({\mathcal U} \coprod {\mathcal V},{\mathcal O}_{\overline X})$. Comme $\overline X$ est donnée par ${\bf P}^1$ avec deux  points $P_{0}$ et $Q_{0}$ identifi\'es, on doit donc avoir $f(P_{0})=f(Q_{0})$ et comme plus haut, on peut montrer qu'il existe un \'element $(h_{0},h_{1}) \in C^{0}({\mathcal U} \coprod {\mathcal V},{\mathcal O}_{{\bf P}^1})$ tel que $\delta (h_{0},h_{1})=f$. Cet \'el\'ement n'est cependant pas forc\'ement dans  $C^{0}({\mathcal U} \coprod {\mathcal V},{\mathcal O}_{\overline X})$ mais l'on a :

$f(P_{0})=h_{0}(P_{0})-h_{1}(P_{0})=h_{0}(Q_{0})-h_{1}(Q_{0})=f(Q_{0})$
donc il existe $\lambda $ tel que
$h_{0}(Q_{0})=h_{0}(P_{0})+\lambda $ et $h_{1}(Q_{0})=h_{1}(P_{0})+\lambda $

Posons $h_{0}'=h_{0}- \lambda (\frac{u-{u_{P_0}}}{{u_{Q_0}}-{u_{P_0}}})$ et  $h_{1}'=h_{1}- \lambda
(\frac{v-v_{P_0}}{v_{Q_0}-v_{P_0}})$. Alors $h_{0}'(Q_{0})=h_{0}'(P_{0})$ et  $h_{1}'(Q_{0})=h_{1}'(P_{0})$ donc ces \'el\'ements sont bien dans $C^{0}({\mathcal U} \coprod {\mathcal V},{\mathcal O}_{\overline X})$ et $h_{0}(u)-h_{1}(u)=f(u)- \lambda (\frac{u-u_{P_0}}{u_{Q_0}-u_{P_0}}-\frac{v-v_{P_0}}{v_{Q_0}-v_{P_0}})$ et donc le coker de $\delta$ est engendr\'e par l'\'el\'ement $(\frac{u-u_{P_0}}{u_{Q_0}-u_{P_0}}-\frac{v-v_{P_0}}{v_{Q_0}-v_{P_0}})$.
\\
Cet \'element est l'image d'un g\'en\'erateur de $H^{1}({\overline X},{\bf C})$, il provient de la boucle dans $\overline X$
form\'ee par l'identification de deux points dans ${\bf P}^1$ (de premier groupe de cohomologie trivial).

Ainsi $W_{0}H^{1}({ X},{\bf C})$ est de rang $1$. On en déduit que $h^{0.0}_{H^{1}( {X},{\bf C})}=1$ et,, la dimension de $H^{1}({ X},{\bf C})$ étant $2$, que $h^{1,1}_{H^{1}({X},{\bf C})}=\text{dim}_{\bf C}\, W_{2}H^{1}({ X},{\bf C})/W_{1}H^{1}({ X},{\bf C})= \text{dim}_{\bf C}\, W_{2}H^{1}({ X},{\bf C})/W_{0}H^{1}({ X},{\bf C})=1$.

Afin de déterminer $\alpha_{1}(X)=\alpha(H^{1}({\overline X},{\bf C}))$, explicitons $F^{1}{\bf H}^{1}({\overline X},[ {\mathcal O}_{\overline X} \stackrel{d}{ \rightarrow } r_{*} {\Omega }_{{\overline X}'}^{1}(\log S)])$. C'est la partie de la cohomologie qui vient du noyau de la flèche :

\hspace{1.5cm}$\delta \,:\,C^{0}({\calu \coprod \calv} ,r_{*}{\Omega}_{{\overline X}'}^{1}(\log S)) \rightarrow C^{1}({\calu \cap \calv},r_{*}{\Omega}_{{\overline X}'}^{1}(\log S))$.\\
 Il suffit donc d'exhiber une forme qui n'est pas exacte $\omega \in C^{0}({\calu \cap \calv} ,r_{*}{\Omega}_{{\overline X}'}^{1}(\log S))$ et de prendre $(\omega,\omega) \in 
C^{0}({\calu \coprod \calv} ,r_{*}{\Omega}_{{\overline X}'}^{1}(\log S))$.
La forme $\omega=(\frac{1}{u-p_{0}}-\frac{1}{u-p_1})du$ g\'en\'ere donc $F^{1}H^{1}(X,{\bf C})$. Déterminons les positions relatives de la filtration de Hodge et sa filtration conjuguée. Comme elle est de niveau $1$, il suffit de trouver la dimension  $t^{1,1}_{H^{1}(X,{\bf C})}=\text{dim}_{\bf C}\,F^{1}H^{1}(X,{\bf C}) \cap {\overline F}^{1}H^{1}(X,{\bf C})$ (qui est $0$ ou $1$). La filtration de Hodge induit une filtration sur le dual du premier groupe de cohomologie par :
\begin{center}
$H^{1}_{DR}(X, {\bf C}) \cong H^{1}_{Betti}(X, {\bf C}) \cong (H_{1}(X, {\bf C}))^{*}=(H_{1}(X, {\bf Z}) \otimes {\bf R})^{*}\otimes_{\bf Z} {\bf C}$
\end{center}
Choisissons une base  de $H_{1}(X, {\bf Z})$ form\'ee de ${\gamma}_0$ et ${\gamma}_1$. Soit ${\gamma}_0$ le
lacet qui serait homologue \`a $0$ dans $X \cup p_0$ et ${\gamma}_1$ un lacet d\'ecrivant la boucle formée par l'identification des deux points. Le resultat ne d\'epends pas du choix des g\'en\'erateurs du $H_1$ puisqu'un choix différent revient à multiplier la matrice des périodes par une matrice à coefficients rationnels et préserve la colinéarité complexe. Int\'egrons donc $\omega $ sur les cycles :
\begin{center}
$<w,{\gamma}_0>=\int_{{\gamma}_0}w=2\pi i$ (c'est le r\'esidu de $\omega $ en $p_1$).\\
$<w,{\gamma}_1>=\int_{{\gamma}_1}w= (\int_{0}^{P_1}-\int_{0}^{Q_1}) (\frac{1}{u-p_{0}}-\frac{1}{u-p_1})du=\lbrack \log
(\frac{u-p_1}{u-p_2}) {\rbrack}^{Q_1}_{P_1}=\log (\frac{Q_{1}-p_1}{Q_{1}-p_2})-\log (\frac{P_{1}-p_1}{P_{1}-p_2})$
\end{center}
(L'int\'egrale entre les points $P_{1}$ et $Q_{1}$ sur la courbe $X'$ d\'esingularis\'ee de $X$ repr\'esente en cohomologie le cycle
cr\'e\'e par l'identification de ces points dans $X'$, boucle non nul-homologue dans $X$).

La quantit\'e $(Q_{1},P_{1},p_{1},p_{2}):=( \frac{Q_{1}-p_{1}}{Q_{1}-p_{2}})/ (\frac{P_{1}-p_{1}}{P_{1}-p_{2}})$ est appel\'ee le
birapport des quatres points $Q_{1},P_{1},p_{1},p_{2}$ ( lorsque au moins trois d'entre eux sont diff\'erents ). Si l'un d'eux est
l'$\infty$ le birapport peut \^etre d\'efini en passant \`a la limite par $(\infty
,P_{1},p_{1},p_{2}):=\frac{P_{1}-p_{2}}{P_{1}-p_{1}}.$

En utilisant la notation précédente on a $<w,{\gamma}_1>=\log (Q_{1},P_{1},p_{1},p_{2})$. Supposons qu'il existe $\lambda
\in {\bf C}$ tel que le vecteur de ${\bf C}^2$ $(<w,{\gamma}_0>,<w,{\gamma}_1>)$ soit proportionnel (complexe) au vecteur obtenu par conjugaison :
\begin{center}
$( 2\pi i,\log (\frac{Q_{1}-p_1}{Q_{1}-p_2})-\log (\frac{P_{1}-p_1}{P_{1}-p_2}))= \lambda \overline{( 2\pi
i,\log (\frac{Q_{1}-p_1}{Q_{1}-p_2})-\log (\frac{P_{1}-p_1}{P_{1}-p_2}))}$
\end{center}
Le birapport est invariant par action de $PGL(1)$ i.e $ \forall \tau \in PGL(1)$, on a l'égalité des birapports\\
$(Q_{1},P_{1},p_{1},p_{2})=(\tau (Q_{1}),\tau (P_{1}),\tau (p_{1}),\tau (p_{2}))$, d'où $\log (Q_{1},P_{1},p_{1},p_{2})=\log(\tau (Q_{1}),\tau (P_{1}),\tau (p_{1}),\tau (p_{2}))$. Comme l'action de $PGL(1)$  sur ${\bf P}^1$ est transitive sur trois points, on peut supposer sans changer la condition de
colin\'earit\'e que l'on a $p_{1}=0, p_{2}=1$ et $P_{1}=\infty$ (on garde la notation $Q_{1}$ pour son image par $\tau$), d'o\`u $\,\,\log(\frac{Q_{1}-p_1}{Q_{1}-p_2})-\log(\frac{P_{1}-p_1}{P_{1}-p_2})=\log(\frac{Q_{1}}{Q_{1}-1})$.
Or la colin\'earit\'e impose que $\lambda =-1$ et donc que $\log(\frac{Q_{1}}{Q_{1}-1})$ soit un imaginaire pur c'est \`a dire que $Q_1$ soit sur la droite $\Re(u)=\frac{1}{2}$. Donc :

$$
\left\lbrace \begin{array}{l}
         t^{1,1}_{H^{1}(X,{\bf C})}=\text{dim}_{\bf C}\,F^{1}H^{1}(X,{\bf C}) \cap {\overline F}^{1}H^{1}(X,{\bf C})=1 \,\,si\,\, Q_1 \in \Re(u)=\frac{1}{2},\\
         t^{1,1}_{H^{1}(X,{\bf C})}=0 \,\, autrement\,\,et \,\, alors\,\,t^{1,0}_{H^{1}(X,{\bf C})}=1 \,\,et\,\,t^{0,1}_{H^{1}(X,{\bf C})}=1.
\end{array}
\right.
$$

Notons ${\mathcal M}_{0,4}$ l'espace des modules des courbes stables de genre $0$ sur ${\bf C}$ avec quatres points marqués. Il est naturellement isomorphe à ${\bf P}^{1}-\{0,1,\infty\}$. Désignons par $X_{Q_1}$ une courbe dans la classe d'isomorphisme représentée par $Q_{1} \in {\mathcal M}_{0,4}$, alors, comme :
\begin{eqnarray*}
\alpha_{1}(X_{Q_1})&=&\frac{1}{2}\Bigl( (0+0)^{2}(h^{0,0}_{H^{1}(X_{Q_1},{\bf C}))}-t^{0,0}_{H^{1}(X_{Q_1},{\bf C}))})+(1+0)^{2}(h^{1,0}_{H^{1}(X_{Q_1},{\bf C}))}-t^{1,0}_{H^{1}(X_{Q_1},{\bf C}))})\\
&+&(0+1)^{2}(h^{0,1}_{H^{1}(X_{Q_1},{\bf C}))}-t^{0,1}_{H^{1}(X_{Q_1},{\bf C}))})+(1+1)^{2}(h^{1,1}_{H^{1}(X_{Q_1},{\bf C}))}-t^{1,1}_{H^{1}(X_{Q_1},{\bf C}))}) \Bigr)
\end{eqnarray*}
Ainsi, le niveau de ${\bf R}$-scindement du premier groupe de cohomologie de $X_{Q_{1}}$ est donné par : 
$$
\left\lbrace \begin{array}{l}
         \alpha_{1}(X_{Q_1})=\alpha(H^{1}(X_{Q_1},{\bf C}))=O \,\, si \,\, Q_{1} \in {\mathcal M}_{0,4}- {\{\Re=\frac{1}{2} \}}\\
         \alpha_{1}(X_{Q_1})=\alpha(H^{1}(X_{Q_1},{\bf C}))=1 \,\, si \,\,Q_{1} \in \{\Re=\frac{1}{2}\}
\end{array}
\right.
$$

\subsubsection{G\'en\'eralisation sur ${\bf P}^1$}

Soit $X={\bf P}^{1}$ avec $m$ points enlev\'es $p_{1},...,p_{m}$ et $2n$ points identifi\'es $P_{1},Q_{1},...,P_{n},Q_{n}$deux à deux (on identifie $P_{k}$ à $Q_{k}$ pour tout $k \in [1,n]$). $H^{1}(S,{\bf C})$ est de rang $m+n-1$, il est engendr\'e par $m+n-1$ \'el\'ements dont $n$ proviennent des boucles form\'ees par l'identification des points et $m-1$ repr\'esentent les lacets autour cr\'e\'es par la non-complétude de $X$. Le calcul de l'hypercohomologie est similaire au précédent, on recouvre ${\overline X}'={\bf P}^{1}$ par les deux ouverts standards $\calu $ $\calv $ munis des coordonnées $u$ et $v$. Comme générateurs de $F^{1}H^{1}(X,{\bf C})$  on peut prendre les \'el\'ements $ \omega_{i}=( \frac{1}{u-p_{i}} -\frac{1}{u-p_{i+1}}) du $ avec $i \in [1,m-1]$. Filtrons l'homologie par la filtration de Hodge :  
D\'esignons par $\gamma_j$ $ j \in [1,m-1]$ les g\'en\'erateurs de $ H_{1}(X,{\bf
C})$ qui sont homologues à zéro dans $X \cup p_j$ et par $\beta_j$, $j \in [1,n]$ les lacets \'el\'ements du $ H_{1}(X,{\bf C})$
repr\'esentants les boucles donn\'ees par $P_{j}=Q_j$ respectivement et ne passant par aucun autre point $P_{k}=Q_k$ pour $k \in
[1,n]$.

Les $\beta_j$ ne sont pas d\'efinis canoniquement mais modulo les $\gamma_j$. Changer la base du premier groupe d'homologie ne modifie pas les conditions de colinéarité dans la matrice des périodes car cela revient à la multiplier par une matrice à coefficients rationnels.\\ 
Filtrons ${\bf C}^{m+n-1}$ par $F^1$, on obtient alors $m-1$ vecteurs :

\hspace{1cm}$(<\alpha_{1}, \omega_{i}>,...,<\alpha_{m-1}, \omega_{i}>,<\beta_{1}, \omega_{i}>,...,<\beta_{n}, \omega_{i}>)$ pour $i \in [1,m-1]$,\\
et on regarde s'ils sont colin\'eaires complexes \`a leurs conjugu\'es. On obtient une matrice $A=(B,C) \in Mat_{(m-1)\times (m+n-1)}({\bf C}) $ telle que les coefficients de $B=(b_{i,j})_{i \in [1,m-1], j \in
[1,m+-1]}$ sont donn\'es par :
$$
\left\lbrace \begin{array}{l}
        b_{i,j}= 2\pi \sqrt{-1}  \,\, si \,\, i=j,\\
        b_{i,j}= -2\pi \sqrt{-1}  \,\, si \,\, i=j+1,\\
        b_{i,j}=0 \,\, \text{sinon}.
\end{array}
\right.
$$
où les coefficients sont donnés par les résidus des formes qui $\omega_i$. Les coefficients de la matrice\\
$C=(c_{i,j})_{i  \in [1,m-1],j \in [m,n+m-1]}$ sont : \hspace{0.2cm} $c_{i,j}=\log (p_{i},p_{i+1},P_{j},Q_{j})$.\\
La matrice des périodes est donc :
$$A=\left(
    \begin{array}{ccccccccc}
     2 \sqrt{-1} \pi  & 0                  & ... & ...              & 0                & \vert   & ... & ... & ... \\
     -2 \sqrt{-1} \pi & 2 \sqrt{-1} \pi    &   0 & ...              & 0                & \vert   & ... & ... & ...           \\
     ...              &...                 &...  &...               &...               & \vert    & ... & ... & ...           \\
      ...             &...                 &...  &...               &...               & \vert     & ... &
c_{ij}=(\log (p_{i},p_{i+1},P_{j},Q_{j})) & ...            \\
      ...             &...                 &...  &...               &...               & \vert      & ... & ... & ...          \\
       0              & ...                & 0   & -2 \sqrt{-1} \pi & 2 \sqrt{-1} \pi  & \vert         & ... & ... & ...
\end{array}
\right)$$
Les vecteurs lignes de la matrice $A$ engendrent l'image de $F^{1}H^{1}(X,{\bf C})$ sur $H_{1}(X,{\bf C})$ et donc calculer $t^{1,1}_{H^{1}(X,{\bf C})}=\emph{dim}_{\bf C}\,F^{1}H^{1}(X,{\bf C}) \cap {\overline F}^{1}H^{1}(X,{\bf C})$ revient à calculer le rang de la matrice $\left(
\begin{array}{c}
A\\
{\overline A}
\end{array}
\right) \in  Mat_{(2m-2)\times (m+n-1)}({\bf C}) $. D'aprés la forme de $B$ ceci revient à voir si chaque vecteur colonne de $A$ est opposé à son conjugué, or :
$$c_{i,j}=-{\overline c}_{i,j} \Leftrightarrow \log (p_{i},p_{i+1},P_{j},Q_{j})=-{\overline {\log (p_{i},p_{i+1},P_{j},Q_{j})}} \Leftrightarrow \vert c_{i,j} \vert =1.$$
Finalement notons ${\mathcal M}_{0,m+2n}$ l'espace des modules des courbes stables de genre $0$ sur ${\bf C}$ avec $m+2n$ points marqués $(p_{1},...,p_{m},P_{1},Q_{1},...,P_{n},Q_{n})$. Il est naturellement isomorphe à $({\bf P}^{1}-\{0,1,\infty\})^{m+n-1} - \Delta$ où $\Delta $ est la grande diagonale définie par l'ensemble des $n-3$-uplets dec points de $({\bf P}^{1}-\{0,1,\infty\})$ tels que deux au moins co\"{\i}ncident. Ainsi :
\begin{prop}Le niveau de ${\bf R}$-scindement de $X$ est donné par :

\hspace{0.7cm}$\alpha_{1}(X)=(m-1)-\text{card} \{ i \in [1,m-1]\, \vert \, \forall j \in [m,m+n-1]\,\,, \vert \log (p_{i},p_{i+1},P_{j},Q_{j}) \vert=1 \}
$
\end{prop}

\textit{Remarques :} $\bullet$ Pour tout $i,j \in [1,m-1] \times [m,m+n-1]$, il existe $\tau \in PGL(1)$ tel que $\tau(p_{i},p_{i+1},P_{j},Q_{j})=(0,1,\infty,\tau(Q_{j}))$ et donc le vecteur ligne $i$ est colinéaire à son conjugué ssi pour tout $j \in  [m,m+n-1]$ il exite un élément $ \tau \in PGL(1)$  tel que $\tau(Q_{j})$ appartienne à la droite de partie réelle $\frac{1}{2}$.\\ 
$\bullet$ Sur les éléments de ${\bf H}^{1}({\overline X},[{\mathcal O}_{\overline X} \stackrel{d}{\rightarrow} r_{*} {\Omega}_{{\overline X}'}^{1}(\log S)])$ qui ne sont pas dans le terme $F^{1}$ de la filtration de Hodge : $F^{1}{\bf H}^{1}({\overline X},[{\mathcal O}_{\overline X} \stackrel{d}{\rightarrow} r_{*} {\Omega}_{{\overline X}'}^{1}(\log S)])$. D'après le travail fait sur ${\bf P}^{1}$ avec deux points identifiés, la partie du $H^{1}(X,{\bf C})$ ne provenant pas des formes logaritmiques est donnée par :\\
\scalebox{0.9}[1]{
$<(\frac{u-u_{P_1}}{u_{Q_1}-u_{P_1}})...(\frac{u-u_{P_n}}{u_{Q_1}-u_{P_n}})-
(\frac{v-v_{P_1}}{v_{Q_1}-v_{P_1}})...(\frac{v-v_{P_n}}{v_{Q_1}-v_{P_n}}),...,(\frac{u-u_{P_1}}{u_{Q_n}-u_{P_1}})...
(\frac{u-u_{P_n}}{u_{Q_n}-u_{P_n}})-
(\frac{v-v_{P_1}}{v_{Q_n}-v_{P_1}})...(\frac{v-v_{P_n}}{v_{Q_n}-v_{P_n}})>$.
}

Ce sont les \'el\'ements de $H^{1}({ X},{\bf C})$ qui proviennent de $H^{1}({\overline X},{\bf C})$, cr\'e\'es par les boucles
sur ${\bf P}^1$ venant de l'identification des points $P_i$ et $Q_i$. Le reste de la cohomologie est donn\'e par la non
compl\`etude de $X$ i.e par les \'el\`ements de $H^{1}({ X'},{\bf C})$, exhibés plus haut.


\subsubsection{Les courbes de genre $1$}

Soit $X$ une courbe alg\'ebrique de type finie sur ${\bf C}$ et de genre arithmétique $1$ avec $m$ points enlev\'es $p_{1},...,p_{m}$ et $2n$ points identifi\'es $P_{1},Q_{1},...,P_{n},Q_{n}$. $\overline X'$ est isomorphe au quotient de ${\bf C}$ muni de la coordonn\'ee $z$ par le r\'eseau ${\Lambda}_{\bf Z}= {\bf Z}+{\bf Z} \tau$ où $\tau \in {\bf C}$ et $Im(\tau )>0$. Cherchons les g\'en\'erateurs de $F^{1}H^{1}(X,{\bf C})$. D'apr\`es [Del3] ce sont les \'el\'ements qui proviennent de $ H^{1}(X',{\bf C})$. 
\begin{prop}\cite{mum}
$\Psi(z)= \sum_{i=1}^{i=k-1} \,\lambda_{i} \, \frac{d}{dz} \log( \theta(z-a_{i}))+cste$ avec $\sum_{i=1}^{i=k-1} \lambda_i=1$ est p\'eriodique pour ${\Lambda}_{\bf Z}$ avec des p\^oles simples aux points $a_i+ \frac{1}{2} (1+ \tau )$ et r\'esidus $\lambda_{i}$, o\`u $\theta $ est la fonction theta sur la courbe elliptique donn\'ee par ${\Lambda}_{\bf Z}$, $\theta(z)= \sum_{n \in {\bf Z}} \exp( \pi \sqrt{-1}
n^{2} + 2 \pi \sqrt{-1} n \tau )$.
\end{prop}
On peut donc prendre comme formes g\'en\'eratrices de $F^{1}H^{1}(S,{\bf C})$ :
$$
\left\lbrace \begin{array}{l}
        \omega_{0}=dz \,\, , \\
        \omega_{i}=d(\log( \theta(z-p_{i}-\frac{1}{2} (1+ \tau )))-\log( \theta(z-p_{i+1}-\frac{1}{2} (1+ \tau )))) \,\, avec \,\, i \in [1,m-1]. 
\end{array}
\right.
$$
Soient $\alpha_{0},\alpha_{1},...,\alpha_{m},\beta_{1},...,\beta_{n}$ des g\'en\'erateurs de $H_{1}(X, {\bf C})$ tels que $\alpha_{j}$ pour $ j \in [1,m]$ soit homologue \`a z\'ero dans $ X \cup p_{j}$ et $\beta_{j}$ est un \'el\'ement representant en homologie la boucle dans $X'$ obtenue par le recollement $ P_{j}=Q_{j}$ (les éléments $\beta_{j}$ sont d\'efinis modulo les $\alpha_{j}$).\\
Fitrons ${\bf C}^{m+n+1}$ par $F^{1}$, cela donne $m$ vecteurs :\\
$(<H_{1},\omega_{i}>)=(<\alpha_{0},\omega_{i}>,<\alpha_{1},\omega_{i}>,...,<\alpha_{m},\omega_{i}>,<\beta_{1},\omega_{i}>,...,<\beta_{n},\omega_{i}>)$ pour $i \in [1,m]$. Les $<\alpha_{j},\omega_{i}>$ sont donn\'es par les r\'esidus des g\'en\'erateurs de $F^{1}H^{1}$ et les $<\beta_{j},\omega_{i}>$ sont donn\'es par int\'egration le long des boucles $\beta_{j}$ i.e sur $X'$, par l'int\'egration de la forme $\omega_{i}$ de $P_j$ \`a $Q_j$ :
\begin{center}
$<\beta_{j},\omega_{i}>= \int_{P_j}^{Q_j} \, \omega_{i}=\int_{P_j}^{Q_j} \,d(\log( \theta(z-p_{i}-\frac{1}{2} (1+ \tau )))-\log(
\theta(z-p_{i+1}-\frac{1}{2} (1+ \tau )))$,\\
i.e $ <\beta_{j},\omega_{i}>=\log( \frac{\theta(Q_{j}-p_{i}-\frac{1}{2} (1+ \tau ))}{\theta(Q_{j}-p_{i+1}-\frac{1}{2}
(1+\tau ))})-\log( \frac{\theta(P_{j}-p_{i}-\frac{1}{2} (1+ \tau ))}{\theta(P_{j}-p_{i+1}-\frac{1}{2} (1+ \tau ))}).$
\end{center}
D'où la matrice $A=(B,C)$ o\`u $B$ est une matrice de dimension $m \times m+1$ et $C$ une matrice $m \times n$.\\
$$A=\left(
    \begin{array}{cccccccccc}
     1  & 0                  & ... & ...              & 0                & \vert   & 0 & ... & 0 \\
     \lambda_{1} & 1    &   -1 & 0             & 0                & \vert   & ... & ... & ...           \\
     ...              & 0                 & 1  & -1               &...               & \vert    & ... & ... & ...           \\
      ...             &...                 &...  &...               &...               & \vert     & ... &
c_{ij}=\log( \frac{\theta(Q_{j}-p_{i}-\frac{1}{2} (1+ \tau ))}{\theta(Q_{j}-p_{i+1}-\frac{1}{2}
(1+\tau ))})-\log( \frac{\theta(P_{j}-p_{i}-\frac{1}{2} (1+ \tau ))}{\theta(P_{j}-p_{i+1}-\frac{1}{2} (1+ \tau ))}) & ...   \\
      ...             &...                 &...  &...               &...               & \vert      & ... & ... & ...          \\
       \lambda_{m-1}              & ...                & 0   & 1  & -1  & \vert         & ... & ... & ...
\end{array}
\right)$$
o\`u $i \in [1,m-1]$ et $j \in [1,n]$ et les $\lambda_{i} $ sont dans ${\bf Z}$.\\

Trouver la valeur de $\alpha_{1}(X)$ revient donc \`a trouver les quadruplets $p_{i},p_{i+1},Q_{j},Q_{j+1}$ tels que :\\
  $\log( \frac{\theta(Q_{j}-p_{i}-\frac{1}{2} (1+ \tau
))}{\theta(Q_{j}-p_{i+1}-\frac{1}{2} (1+\tau ))})-\log( \frac{\theta(P_{j}-p_{i}-\frac{1}{2} (1+ \tau
))}{\theta(P_{j}-p_{i+1}-\frac{1}{2} (1+ \tau ))})= \overline{\log( \frac{\theta(Q_{j}-p_{i}-\frac{1}{2} (1+ \tau
))}{\theta(Q_{j}-p_{i+1}-\frac{1}{2} (1+\tau ))})-\log( \frac{\theta(P_{j}-p_{i}-\frac{1}{2} (1+ \tau
))}{\theta(P_{j}-p_{i+1}-\frac{1}{2} (1+ \tau ))})}$ suivant la configuration des points $p_{i},p_{i+1},Q_{j},Q_{j+1}$. D'où :

\begin{prop} Le niveau de ${\bf R}$-scindement de $X$ est donné par :\\
$\alpha_{1}(X)=(m-1)-\text{card} \{ i \in [1,m-1]\, \vert \, \forall j \in [1,n]\,\,, \vert \log( \frac{\theta(Q_{j}-p_{i}-\frac{1}{2} (1+ \tau
))}{\theta(Q_{j}-p_{i+1}-\frac{1}{2} (1+\tau ))})-\log( \frac{\theta(P_{j}-p_{i}-\frac{1}{2} (1+ \tau
))}{\theta(P_{j}-p_{i+1}-\frac{1}{2} (1+ \tau ))})= \overline{\log( \frac{\theta(Q_{j}-p_{i}-\frac{1}{2} (1+ \tau
))}{\theta(Q_{j}-p_{i+1}-\frac{1}{2} (1+\tau ))})-\log( \frac{\theta(P_{j}-p_{i}-\frac{1}{2} (1+ \tau
))}{\theta(P_{j}-p_{i+1}-\frac{1}{2} (1+ \tau ))})} \}
$.

\end{prop}


\textsc{o. penacchio: Laboratoire Emile Picard, Universit\'e Paul Sabatier,
118 route de Narbonne, 31062 Toulouse Cedex France.} \\
E-mail: penacchi@picard.ups-tlse.fr


\begin{thebibliography}{papouas}


\bibitem[C]{car} \textsc{Carlson}, \emph{ Extensions of Mixed Twistor Structures}, Journ\'ees de G\'eom\'etrie Alg\'ebrique d'Angers 1979, A.Beauville (Editor), Sijthoff \& Noordhoff, Alphen aan den Rijn, 1980, pp. 107-128.

\bibitem[C-K-S]{catkapsch} \textsc{Cattani E.,Kaplan A.,Schmid W.}:\emph{ Degeneration of Hodge Structures } Ann.of Math. 123.(1986),457-535.

\bibitem[Del2]{del2} \textsc{Deligne}, \emph{ Th\'eorie de Hodge II } Publ. Math. IHES 40 (1972), 5-57.

 \bibitem[Del3]{del3} \textsc{Deligne}, \emph{ Th\'eorie de Hodge III } Publ. Math. IHES 40 (1975), 6-77.

\bibitem[G-H]{grihar} \textsc{Griffiths P. Harris J.} \emph{ Principles of Algebraic Geometric } Wiley, New York.



\bibitem[H]{har} \textsc{Hartshorne}, \emph{ Algebraic Geometry} Springer GTM 52,1977.


\bibitem[M]{mor} \textsc{Morgan}, \emph{The Algebraic Topology of smooth Algebraic Varieties} Publ. Math. I.H.E.S., 48 (1978), 137-204.


 \bibitem[Mum]{mum} \textsc{Mumford}, \emph{ Tata Lectures on Theta I }, Progress in Mathematics, Birkh\"auser.



        \bibitem[OSS]{oss}
      \textsc{Okonek, Schneider, Spindler},
     \emph{Vector Bundles on Complex Projective Spaces.}, Progress in Mathematics, Vol. 3, Birkha\"user.
 
\bibitem[P]{pen} \textsc{O.Penacchio}, \emph{ Niveau de ${\bf R}$-scindement et variations de structures de Hodge mixtes}, en préparation.

\bibitem[Si]{sim} \textsc{C.T.Simpson}, \emph{ Mixed twistor structures}, math/alg-geom/9705006 .

   


  










































































































 











 \end{thebibliography}
\end{document}